\documentclass[reqno,fleqn,11pt]{amsart}

\usepackage[english]{babel}

\usepackage{amsmath}
\usepackage{amssymb}
\usepackage{amsfonts}
\usepackage{amsthm}

\newcommand{\N}{\mathbb{N}}
\newcommand{\R}{\mathbb{R}}

\newcommand{\sphere}{\mathbb{S}}
\newcommand{\f}{\rightarrow}

\newcommand{\de}{\partial}
\newcommand{\pa}{\parallel}

\newcommand{\defeq}{\mathrel{\mathop:}=}

\newcommand{\jac}{\operatorname{Jac}}
\newcommand{\inj}{\operatorname{inj}}

\newcommand{\Ball}{\operatorname{B}}
\newcommand{\Sob}{\operatorname{Sob}}

\newcommand{\tr}{\operatorname{tr}}

\newcommand{\adeg}{\operatorname{Adeg}}
\newcommand{\gdeg}{\mathcal G_{\mathrm{deg}}}

\newcommand{\Id}{\operatorname{Id}}

\newcommand{\Vol}{\operatorname{Vol}}

\newcommand{\Ric}{\operatorname{Ricci}}
\newcommand{\diam}{\operatorname{Diam}}
\newcommand{\can}{\operatorname{can}}

\newtheorem{thm}{Theorem}[section]

\newtheorem{prop}[thm]{Proposition}
\newtheorem{propA1}{Proposition A.}
\newtheorem{lemma}[thm]{Lemma}
\newtheorem{lemmaA}[propA1]{Lemma A.}
\newtheorem{lemmaB}{Lemma B.}
\newtheorem{lemmaC}{Lemma C.}

\theoremstyle{remark}
\newtheorem{exmp}[thm]{Example}
\newtheorem{rmk}[thm]{Remark}

\theoremstyle{definition}
\newtheorem{defn}[thm]{Definition}

\begin{document}

\title{A  spectra comparison theorem and its applications}
\author[F.Cerocchi]{Filippo Cerocchi}
\address{Dipartimento di Matematica "G.Castelnuovo", Universit\`{a} di Roma "Sapienza",
Piazzale Aldo Moro 5, 00185 Roma, Italy}
\address{Institut Fourier, Universit\'e Grenoble 1, 100 rue des maths, BP 74, 38402 St. Martin d'H\`eres, France\newline}

\email{\newline cerocchi@mat.uniroma1.it;   filippo.cerocchi@ujf-grenoble.fr; }

\begin{abstract}
We give a sharp comparison between the spectra of two Riemannian manifolds $(Y,g)$ and $(X,g_0)$ under the following assumptions: $(X,g_0)$ has bounded geometry, $(Y,g)$ admits a continuous Gromov-Hausdorff $\varepsilon$-approximation onto $(X,g_0)$ of non zero absolute degree, and the volume of $(Y,g)$ is almost smaller than the volume of $(X,g_0)$. These assumptions imply no restrictions on the local topology or geometry of $(Y,g)$ in particular no curvature assumption is supposed or infered.
\end{abstract}

\maketitle
 
\section{Introduction}

The aim of this paper is to compare the spectra of two Riemannian manifolds $(Y,g)$ and $(X,g_0)$ and to bound the gap in terms of the Gromov-Hausdorff distance between these two spaces when this distance is smaller than some universal constant (see Theorem \ref{SpC-MT2}).\\\\
Estimates from above and from below for the eigenvalues of the Laplace-Beltrami operator of manifolds satisfying a lower bound of the Ricci curvature and an upper bound of the diameter were derived in the decade from 1975 to 1985. Namely, following  S. Y. Cheng \cite{Cheng1}, and P. Li and S. T. Yau, \cite{Li-Yau}, we know that, when $(Y,g)$ is a compact Riemannian $n$-manifold of diameter $\diam(Y,g)\leq D$, whose Ricci curvature, $\Ric_g$, satisfies the bound $\Ric_g\geq-(n-1)\kappa^2$, then the eigenvalues of the Laplace-Beltrami operator admit the following upper bound:
\begin{equation}\label{SpC-Ubound}
\lambda_k(Y,g)\leq\frac{C(\alpha)}{\Vol_g(Y)^{\frac{2}{n}}}\cdot k^{\frac{2}{n}}
\end{equation}
where $\alpha=\kappa\,D$. 
On the other hand under the same assumptions we have the following inferior bound for the eigenvalues of the Laplace-Beltrami operator of $(Y,g)$ (see P. Li and S. T. Yau \cite{Li-Yau}, M. Gromov \cite{Gromov} and S. Gallot, \cite{Ga1}, \cite{Ga2}, \cite{Ga3}):
\begin{equation}\label{SpC-Lbound}
\lambda_k(Y,g)\geq\frac{\Gamma(\alpha)}{\diam(Y,g)^2}\cdot k^{\frac{2}{n}}.
\end{equation}
Explicit values for the constants $C(\alpha)$ and $\Gamma(\alpha)$ can be found in \cite{Li-Yau}, \cite{Ga1}, \cite{Ga2}, \cite{Ga3} or \cite{BBG}.

\begin{rmk}
As long as we are only concerned by their dependence with respect to the index $k$, the inequalities (\ref{SpC-Ubound}) and (\ref{SpC-Lbound}) agree with the well known Weyl asymptotic formula which says that the sequence of the eigenvalues of the Laplace-Beltrami operator of a compact Riemannian $n$-manifold $(Y,g)$ behaves asymptotically like 
$\lambda_k\sim\frac{(2\pi)^2}{(\Vol(\mathrm{B}^n)\cdot\Vol_g(M))^{\frac{2}{n}}}\cdot k^{\frac{2}{n}}\,.$
\end{rmk}

For the same reason, assume that we have two compact Riemannian $n$-manifolds $(Y,g)$, $(X,g_0)$, whose Ricci curvatures satisfy the same lower bound $$\Ric_g\geq-(n-1)\kappa^2,\;\Ric_{g_0}\geq-(n-1)\kappa^2$$
and such that $\diam(Y,g), \diam(X,g_0)\leq D$. From inequalities (\ref{SpC-Ubound}) and (\ref{SpC-Lbound}) we obtain:
$$\lambda_k(Y,g)\leq\frac{C(\alpha)}{\Gamma(\alpha)}\cdot\left(\frac{\diam(X,g_0)^n}{\Vol_g(Y)}\right)^{\frac{2}{n}}\cdot\lambda_k(X,g_0)$$
and, exchanging the roles of $(Y, g)$ and $(X,g_0)$:
$$\lambda_k(X,g_0)\leq\frac{C(\alpha)}{\Gamma(\alpha)}\cdot\left(\frac{\diam(Y,g)^n}{\Vol_{g_0}(X)}\right)^{\frac{2}{n}}\cdot\lambda_k(Y,g)$$
However, since the quantity $C(\alpha)/\Gamma(\alpha)$ is considerably greater than $1$ (although the constants involved are sharp!) and as the ratio $\frac{(\diam)^n}{\Vol}$ can be arbitrarily large, even if we suppose that the diameters and the volumes of $(Y,g)$ and $(X,g_0)$ are almost the same, we cannot infer an  equality  or deduce some sharp pinching result between the $k^{th}$ eigenvalues.\\

In \cite{ChCo}, Theorem 7.11,  J. Cheeger and T. Colding gave a convergence result for the eigenvalues of the Laplace operators of a sequence of $n$-dimensional manifolds $(Y_k,g_k)$ whose Ricci curvatures are bounded from below by $-(n-1)$ and which converge with respect to the Gromov-Hausdorff distance to a given smooth manifold $(X,g_0)$ of the same dimension (notice that, $(X,g_0)$ being fixed, its Ricci curvature is automatically bounded from below). Namely, they prove that, for any fixed $j\in\N^*$, the $j^{th}$ eigenvalue $\lambda_j(Y_k,g_k)$ of the Laplace-Beltrami operator of $(Y_k,g_k)$ converges to $\lambda_j(X,g_0)$ as $k\f+\infty$.\\
Notice that the Gromov-Hausdorff convergence is not, in itself, a strong assumption (in particular it gives no informations on the local geometries of the $(Y_k,g_k)$), but that it becomes quite a strong one when it is combined with a uniform lower bound on the Ricci curvature of the Riemannian manifolds $(Y_k,g_k)$. Assuming together these two properties one obtains that, for large values of $k$, the local geometries of $(Y_k,g_k)$ are almost the same (J. Cheeger, T. Colding \cite{ChCo} and T. Colding \cite{Col}), in particular, for every $\varepsilon>0$, if $y_k\f x$, the volume of the geodesic ball $B(y_k,\varepsilon)$ of $(Y_k,g_k)$ converges to the volume of the geodesic ball $B(x,\varepsilon)$ of $(X,g_0)$. Moreover $Y_k$ is diffeomorphic to $(X,g_0)$ for large values of $k$ (\cite{ChCo}), and $\Vol_{g_k}(Y_k)\f\Vol_{g_0}(X)$ as $k\f+\infty$ (\cite{Col}). The fact that $\diam(Y_k,g_k)$ converges to $\diam(X,g_0)$ is an immediate consequence of the Gromov-Hausdorff convergence.\\

We remark that the eigenvalues approximation methods show that we have convergence of $\lambda_i(Y_k,g_k)\f\lambda_i(X,g_0)$ when $(Y_k,g_k)$ is a sequence of polyedral approximations converging to $(X,g_0)$ (see \cite{dodpat}, \S3). However the result of Dodziuk and Patodi does not provide an upper bound of the "error" $|\lambda_i(Y_k,g_k)-\lambda_i(X,g_0)|$.\\

Comparing with the above mentioned results, the comparison between the spectra of two manifolds $(Y,g)$ and $(X,g_0)$ that we aim must obey to quite different rules: namely we are authorized to assume that the geometry of $(X,g_0)$ is bounded. On the contrary, on $(Y,g)$, any assumption which implies a control on the local topology or geometry is prohibited. Let us denote by $\sigma_0$ the sectional curvature of $(X,g_0)$ and by $\inj(X,g_0)$ its injectivity radius,  the main result in this direction is the following:

\begin{thm}\label{SpC-MT2}
Let $(X^n,g_0)$ be a compact, connected, Riemannian manifold satisfying the assumptions:
$$\diam(X,g_0)\leq D,\quad\inj(X,g_0)\geq i_0,\quad|\sigma_0|\leq\kappa^2,$$
where $D,\, i_0,\,\kappa$ are arbitrary positive constants; let $\varepsilon_0=\varepsilon_0(n,i_0,\kappa)$ be the universal constant given in Theorem \ref{SpC-barycenter}.\\
Let $(Y^n,g)$ be any compact, connected Riemannian manifold such that there exists a continuous Gromov-Hausdorff $\varepsilon$-approximation $f:(Y,g)\f(X,g_0)$ of non zero absolute degree, where
\begin{equation}\label{SpC-*}
\varepsilon<\varepsilon_1(n,i_0,\kappa)=\min\left\{\varepsilon_0(n,i_0,\kappa);\,\frac{1}{\kappa}\left(\frac{\left(\frac{10}{9}\right)^{\frac{2}{n}}-1}{20(n+1)}\right)^4\right\}
\end{equation}
If we assume that
$$
[1-10\,n\,(n+1)\,(\kappa\varepsilon)^{\frac{1}{4}}]\cdot\Vol_g(Y)<\Vol_{g_0}(X)
$$
then, for every $i\in\N$, we have
\begin{equation*}
\lambda_i(Y,g)\leq\left(1+C_1(n)(\kappa\varepsilon)^{\frac{1}{16}}\right)\cdot
\end{equation*}
\begin{equation}\label{SpC-2*}
\cdot\left(1+C_2(n,\,\kappa\, D,\,D^2\cdot\lambda_i(X,g_0))(\kappa\varepsilon)^{\frac{1}{8}}\right)\cdot\lambda_i(X,g_0)
\end{equation}
where
$$C_1(n)=14\,(n-1)\cdot\sqrt{n+1}$$
$$C_2(n,\,\alpha,\Lambda)=4(n+1)\,\left[(2n+1)\,e^n\,\left[1+B(\alpha)\sqrt{\Lambda+(n-1)\alpha^2}\right]^n+2\right]$$
where $B(\alpha)$ is the isoperimetric constant defined in Proposition \ref{SpC-GSobolev} and where the right hand side of (\ref{SpC-*}) goes to $\lambda_i(X,g_0)$ when $\varepsilon\f0_+$.
\end{thm}

\begin{rmk}
\begin{enumerate}
\item The inequality (\ref{SpC-*}) given by Theorem \ref{SpC-MT2} is sharp: in fact it provides an upper bound of $\frac{\lambda_i(Y,g)}{\lambda_i(X,g_0)}$ which goes to $1$ as $\varepsilon\f 0^+$.\\
\item Notice that the only assumptions that we make on $(Y,g)$ in the Theorem \ref{SpC-MT2} are:\\
\begin{itemize}
\item[(i)] the existence of a continuous Gromov-Hausdorff approximation of nonzero absolute degree from $(Y,g)$ to $(X,g_0)$;\\
\item[(ii)] the assumption that the volume of $(Y,g)$ is almost smaller than the volume of $(X,g_0)$.\\
\end{itemize}
Notice that there is no assumption on the curvature of $(Y,g)$.\\
The weakness of these assumptions on $(Y,g)$ is first illustrated by the fact that they give no information on the local topology of $(Y,g)$ or on the local topology of $Y$. In fact, in the Example \ref{SpC-4.2}, for every $(X^n,g_0)$ we construct a family of pairwise non homotopic Riemannian manifolds $(Y_\varepsilon,g_\varepsilon)$ which satisfy the assumptions (i) and (ii) above (with $\varepsilon\f0_+$) and thus converge to $(X,g_0)$ as $\varepsilon\f0_+$.\\
\item Another illustration of the weakness of the assumptions made on $(Y,g)$ is the fact that it is impossible to get a lower bound of $\frac{\lambda_i(Y,g)}{\lambda_i(X,g_0)}$ under these assumptions. In fact, in the Example \ref{SpC-4.3} we construct, for any fixed Riemannian manifold $(X,g_0)$, a sequence of Riemannian manifolds $(Y_k,g_k)$ (diffeomorphic to $(X,g_0)$), which satisfy assumptions (i) and (ii) above, and such that $\frac{\lambda_1(Y_k,g_k)}{\lambda_1(X,g_0)}\f0$ when $k\f+\infty$.\\
Let us stress the fact that the counter-examples mentioned above in the remarks (3) and (4) satisfy all the assumptions of the Theorem \ref{SpC-MT2} for arbitrarily small values of $\varepsilon$.\\
\item The Theorem \ref{SpC-MT2} is not only a convergence result: in fact, it is valid for non small values of $\varepsilon$ (\textit{i.e.} for every $\varepsilon<\varepsilon_1(n,i_0,\kappa)$). For every $\varepsilon<\varepsilon_1(n,i_0,\kappa)$, it provides an explicit upper bound for the "error term" $\left(\frac{\lambda_i(Y,g)-\lambda_i(X,g_0)}{\lambda_i(X,g_0)}\right)$.\\
\item Theorem \ref{SpC-MT2} also works when $g$ is not a smooth Riemannian metric (for example if $g$ is piecewise $C^1$). It thus provides a sharp estimate of $\lambda_i(X,g_0)$ by the corresponding eigenvalue of a polyedral $\varepsilon$-approximation and a bound of the error (in one sense).
\end{enumerate}
\end{rmk}

Theorem \ref{SpC-MT2}  is a consequence of Theorem \ref{SpC-barycenter} and of the following technical result:
\begin{prop}\label{SpC-tech1} Let $\kappa, D >0$.
Let $(X,g_0)$ be a connected, compact Riemannian manifold which satisfy $\Ric_{g_0}\geq-(n-1)\kappa^2$ and $\diam(X,g_0)\leq D$. Let $(Y,g)$ be another compact, connected Riemannian manifold such that $(1-\eta)\Vol_g(Y)\leq\Vol_{g_0}(X)$ (where $0<\eta\leq\frac{1}{9}$) and such that there exists a Lipschitz map $F:(Y,g)\f(X,g_0)$ of non zero absolute degree which verifies the following bound on the pointwise energy:
$e_y(F)\leq n(1+\eta)^{2/n}$ \textit{a.e.}.
Then
$$\lambda_i(Y,g)\leq(1+7(n-1)\,\eta^{\frac{1}{4}})\cdot(1+C(n, D^2\lambda_i(X,g_0),\alpha)\,\eta^{\frac{1}{2}})\cdot\lambda_i(X,g_0)$$
where $C(n,D^2\lambda_i,\alpha)=(2n+1)\cdot e^{n}\cdot\left(1+B(\alpha)\sqrt{\lambda_i D^2+(n-1)\alpha^2}\right)^n+2$ and where $\alpha=\kappa\cdot D$.
\end{prop}

\begin{rmk}[Dimension $n=2$]
In order to simplify the notations we shall give the proofs only for dimensions $n>2$. However, the same arguments hold in dimension $n=2$, provided some slight modifications. Just observe that:
\begin{itemize}
\item[(a)] Lemma \ref{SpC-A}, (i)  is valid for $n=2$;
\item[(b)] Lemma \ref{SpC-A} (ii) is  valid if we replace $n$ by $p$ where $p=n$ for $n>2$ and $p>n$ for $n=2$ (it is sufficient to replace, in the definition of $h$, the function $f^{\frac{2(n-1)}{n-2}}$ by $f^{\frac{2(p-1)}{p-2}}$);
\end{itemize}
all the arguments then works, substituting $n$ by $p$, included the Sobolev inequality (which is not sharp in  dimension $2$) which says that there exists a constant $C$ such that
$$\pa f\pa_{\frac{2p}{p-2}}\leq C\cdot\pa df\pa_{2}+\pa f\pa_2$$
Moser's iteration method then works with $\beta=\frac{p}{p-2}$. 
\end{rmk}

\section{Geometric-analytic tools}

Let us recall the following definition:

\begin{defn} Let $(M,g)$ be a Riemannian manifold. Let us denote by $\Ric_g$ the Ricci curvature of $(M,g)$.
We define the invariant $r_{\min}$ as the infimum of $\Ric_g$ viewed as function on the unit tangent bundle $U_gM$.
\end{defn}

We shall use a Sobolev inequality due to S. Gallot. The original results of S. Gallot were published in a short note in the \emph{Comptes rendus de l'Acad\'emie de Sciences} (see \cite{Ga1}); the original proofs are rather dense, thus, for easier reference, we shall present the detailed proofs in Appendix.

\begin{prop}[Sobolev inequality, \cite{Ga1}]\label{SpC-GSobolev}
Let $(M,g)$ be a compact Riemannian manifold of dimension $n$, such that $\diam(M,g)\leq D$,  $r_{\min}\cdot D^2\geq-(n-1)\alpha^2$. For every function $f:M\f\R$ in $H_1^2(M,g)$ we have:
\begin{equation*}
\mathrm{(i)}\quad\left(\frac{1}{V}\int_M|f-\overline f|^{\frac{2n}{n-2}}\;dv_g\right)^{\frac{n-2}{2n}}\leq
\end{equation*}
$$\leq\left[\frac{2(n-1)}{(n-2)\Gamma(\alpha)}+\frac{2}{H(\alpha)}\right] \cdot D\cdot\left(\frac{1}{V}\int_M|\nabla f|^2\,dv_g\right)^{\frac{1}{2}}$$
\begin{equation*}
\mathrm{(ii)}\quad\left(\frac{1}{V}\int_M f^{\frac{2n}{n-2}}\,dv_g\right)^{\frac{n-2}{2n}}\leq
\end{equation*}
$$\leq\left(\frac{2(n-1)}{(n-2)\Gamma(\alpha)}+\frac{2}{H(\alpha)}\right)\cdot D\cdot\left(\frac{1}{V}\int_M|\nabla f|^2\,dv_g\right)^{\frac{1}{2}}+\left(\frac{1}{V}\int_M f^2\,dv_g\right)^{\frac{1}{2}}$$
where we denote by $\bar f$ the mean value of $f$, \textit{i.e.} $\bar f=\frac{1}{V}\int_M f\,dv_g$, where  $H(\alpha)=\alpha\left(\int_0^{\alpha/2}(\cosh(t))^{n-1}dt\right)^{-1}$ and  where 
$$\Gamma(\alpha)=\alpha\left(\int_0^\alpha\left(\frac{\alpha}{H(\alpha)}\cosh(t)+\frac{1}{n}\sinh(t)\right)^{n-1}\,dt\right)^{-\frac{1}{n}}$$ We will use the notation $B(\alpha)$ to refer to the quantity $\left(\frac{2(n-1)}{(n-2)\Gamma(\alpha)}+\frac{2}{H(\alpha)}\right)$.
\end{prop}

\subsection{Moser's iteration method}
The method that we are going to use in this section has been introduced by J. Moser in 1961 (see \cite{Mos}) in order to prove a Harnack's inequality for solutions of second order, uniformly elliptic partial differential equation of selfadjoint form. The original method allows to derive $L^\infty$ estimates, for eigenfunctions of a differential operator of the prescribed type, in terms of the geometric data of the domain under consideration, this is achieved by means of an iterated use of a Sobolev inequality (see \cite{Mos}, \S 4).

Moser's iteration method has been widely used in Spectral Geometry to obtain eigenvalues estimates for the Laplace-Beltrami operator, the Hodge-de Rham Laplacian and the $p$-Laplacian under appropriate geometric assumptions, see \cite{Li}, \cite{Ga1}, \cite{Ga2} and, for a more recent application, \cite{ACGR}.\newline

We shall use the same symbol $\Delta_g$ to denote on one hand the usual Laplace-Beltrami operator on functions and, on the other hand, the Hodge-de Rham Laplacian viewed as an operator acting on the space $d[C^\infty(M)]$ of exact differential forms of degree $1$ (the discrimination between these two cases will be given by the context). As $\Delta_g$ commutes with the exterior derivative $d$, it comes that the Hodge-de Rham Laplacian maps $d[C^\infty(M)]$ onto $d[C^\infty(M)]$ and that $d$ maps each eigenspace of the Laplace operator on the eigenspace of the Hodge-de Rham Laplacian (restricted to $d[C^\infty(M)]$) corresponding to the same eigenvalue.\\ Let us define $\mathcal A(\lambda)$ as the direct sum of the eigenspaces of the Laplace-Beltrami operator (acting on functions) corresponding to the eigenvalues $\lambda_i\leq\lambda$, the above commutation implies that $\mathcal A_1(\lambda) \defeq d[\mathcal A(\lambda)]$ is also a direct sum of the eigenspaces of the Hodge-de Rham Laplacian (acting on $d[C^\infty(M)]$) corresponding to the eigenvalues $\lambda_i\leq\lambda$.\\
 In order to obtain a clearer statement for the next proposition, let $\beta=\frac{n}{n-2}$ and let us define the function:
$$\xi(x)=\prod_{i=0}^{\infty}\left(1
+\frac{\beta^i}{\sqrt{2\beta^i-1}}\cdot x\right)^{\beta^{-i}}$$

\begin{rmk}
We remark here that the infinite product defining $\xi$ is convergent, as proved in Appendix B.
\end{rmk}

\begin{prop}[revisiting \cite{Ga2}]\label{SpC-GSobolev2}
Let $(M,g)$ be a compact Riemannian manifold such that $$r_{\min}\cdot\diam(M,g)^2\geq-(n-1)\alpha^2$$ and such that $\diam(M,g)\leq D$. For any function $f\in\mathcal A(\lambda)$, we have:
\begin{equation*}
\mathrm{(i)} \pa f\pa_\infty\leq\xi(B(\alpha)D\sqrt\lambda)\cdot\pa f\pa_2\leq
\end{equation*}
$$
\leq\exp\left(\frac{n}{2}\cdot\frac{B(\alpha)D\sqrt\lambda}{1+B(\alpha)D\sqrt\lambda}\right)\cdot\left(1+B(\alpha)D\sqrt\lambda\right)^{\frac{n}{2}}\cdot\pa f\pa_2$$
\begin{equation*}
\mathrm{(ii)} \pa df\pa_\infty\leq\xi\left(B(\alpha)\sqrt{\lambda D^2+(n-1)\alpha^2}\right)\cdot\pa df\pa_2\leq
\end{equation*}
$$\leq\exp\left(\frac{n}{2}\cdot\frac{B(\alpha)\sqrt{\lambda D^2+(n-1)\alpha^2}}{1+B(\alpha)\sqrt{\lambda D^2+(n-1)\alpha^2}}\right)\cdot$$$$\cdot\left(1+B(\alpha)\sqrt{\lambda D^2+(n-1)\alpha^2}\right)^{\frac{n}{2}}\cdot\pa df\pa_2$$
where $B(\alpha)$ is the Sobolev constant that we defined in the statement of Proposition \ref{SpC-GSobolev}.
\end{prop}

\begin{rmk}\label{SpC-rmkGSobolev2}
The quantities $\alpha$ and $\lambda D^2$ in the previous statement are invariant under homotheties.
\end{rmk}

\textbf{Proof of Proposition \ref{SpC-GSobolev2}.} Let $S$ denote any element of $\mathcal A_1(\lambda)$ (resp. of $\mathcal A(\lambda)$), such that $S=df$ (resp. $S=f$) for some $C^\infty$ function $f$. As we have supposed that $r_{\min}\geq-(n-1)\kappa^2$ (where $\kappa=\alpha/\diam(M,g)$), we may introduce a new constant $\kappa_0$, which allows to handle both cases (\textit{i.e.} the case where $S\in\mathcal A_1(\lambda)$ and the case where $S\in\mathcal A(\lambda)$) in a unique computation: we thus define $\kappa_0$ by:
$$\kappa_0=\left\{
\begin{array}{c}
\kappa \mbox{ when }\Delta_g\mbox{ is the Hodge-de Rham Laplacian and when } S\in\mathcal A_1(\lambda);\\
0 \mbox{ when }\Delta_g\mbox{ is the Laplace-Beltrami operator and when }S\in\mathcal A(\lambda).
\end{array}\right.
$$\\
We precise that we use the notation $|S|$ or $|S|(x)$ to denote the pointwise norm of $S$, whereas we use the notation $\pa S\pa_p$ when we consider the global $L^p$ norm of $S$.\\
By the B\"ochner Formula (resp. by definition of the Laplace-Beltrami operator) we have:
$$\langle\nabla^*\nabla S,S\rangle=\langle\Delta_g S, S\rangle -\Ric_g(\nabla f,\nabla f)\leq |\Delta_g S|\cdot|S|+(n-1)\,\kappa^2\,|S|^2$$
when $S\in\mathcal A_1(\lambda)$ (resp. $\langle\nabla^*\nabla S,S\rangle=\Delta_g S\cdot S$, when $S\in\mathcal A(\lambda)$), which gives the following formula:
\begin{equation}\label{SpC-GS2eq1}
\langle\nabla^*\nabla S,S\rangle\leq\left(|\Delta_g S|+(n-1)\,\kappa_0^2\,|S|\right)|S|
\end{equation}
which is valid in both cases (the case $S\in\mathcal A_1(\lambda)$ and the case $S\in\mathcal A(\lambda)$) by definition of $\kappa_0$.
We define (in the case $S\in\mathcal A_1(\lambda)$ as in  the case $S\in\mathcal A(\lambda)$), the function $s_\varepsilon$ by $s_\varepsilon=\sqrt{|S|^2+\varepsilon^2}$. A direct computation leads to 
$|ds_\varepsilon|^2\leq|\nabla S|^2$ and $\Delta_g (s_\varepsilon^2)=\Delta_g(|S|^2)$, so that:
$$s_\varepsilon\,\Delta_g s_\varepsilon=\frac{1}{2}\Delta_g(s_\varepsilon^2)+|d s_\varepsilon|^2\leq\frac{1}{2}\Delta_g(|S|^2)+|\nabla S|^2=\langle\nabla^*\nabla S,S\rangle\leq$$
$$\leq(|\Delta_g S|+(n-1)\,\kappa_0^2\,|S|)|S|\leq(|\Delta_g S|+(n-1)\,\kappa_0^2\,|S|)s_\varepsilon$$
For every $k>\frac{1}{2}$, we deduce that
$$\int_M|d(s_\varepsilon^k)|^2\,dv_g=\frac{k^2}{2k-1}\int_M\langle ds_\varepsilon, d(s_\varepsilon^{2k-1})\rangle\, dv_g=$$
$$=\frac{k^2}{2k-1}\int_M s_\varepsilon^{2k-1}(\Delta_g s_\varepsilon)\,dv_g\leq$$
$$\leq\frac{k^2}{2k-1}\left[\int_M|\Delta_g S|\cdot s_\varepsilon^{2k-1}\,dv_g+(n-1)\,\kappa_0^2\,\int_M|S|\cdot s_\varepsilon^{2k-1}\,dv_g \right]\leq$$
$$\leq\frac{k^2 V}{2k-1}\left[\pa\Delta_g S\pa_{2k}\cdot\pa s_\varepsilon\pa_{2k}^{2k-1}+(n-1)\,\kappa_0^2\,\pa S\pa_{2k}\cdot\pa s_\varepsilon\pa_{2k}^{2k-1}\right]$$
It follows that:
\begin{equation*}
\left(\frac{1}{V}\int_M|d(s_\varepsilon^k)|^2\,dv_g\right)^{\frac{1}{2}}\leq
\end{equation*}
\begin{equation}\label{SpC-GS2eq2}
\leq\frac{k}{\sqrt{2k-1}}\pa s_\varepsilon\pa_{2k}^{k-\frac{1}{2}}\cdot\left(\pa\Delta_g S\pa_{2k}+(n-1)\kappa_0^2\pa S\pa_{2k}\right)^{\frac{1}{2}}
\end{equation}
On the other hand the Sobolev inequality of Proposition \ref{SpC-GSobolev} (ii) gives:
$$\left(\frac{1}{V}\int_M|d(s_\varepsilon^k)|^2\,dv_g\right)^{\frac{1}{2}}\geq\frac{1}{B(\alpha) D}\left[\pa s_\varepsilon^k\pa_{\frac{2n}{n-2}}-\pa s_\varepsilon^k\pa_2\right]$$
hence, putting together this estimate with inequality (\ref{SpC-GS2eq2}), we get:
\begin{equation*}
\pa s_\varepsilon\pa_{\frac{2kn}{n-2}}^k-\pa s_\varepsilon\pa_{2k}^k\leq
\end{equation*}
\begin{equation}\label{SpC-GS2eq3}
\leq\frac{[B(\alpha)D]\cdot k}{\sqrt{2k-1}}\cdot\pa s_\varepsilon\pa_{2k}^{k-\frac{1}{2}}\cdot\left[\pa\Delta_g S\pa_{2k}+(n-1)\kappa_0^2\pa S\pa_{2k}\right]^{\frac{1}{2}}
\end{equation}
Observe that $\pa s_\varepsilon\pa_{2p}\f\pa S\pa_{2p}$ when $\varepsilon\f 0_+$. Thus inequality (\ref{SpC-GS2eq3}) becomes:
\begin{equation*}
\pa S\pa_{\frac{2kn}{n-2}}^k\leq\frac{[B(\alpha)D]\cdot k}{\sqrt{2k-1}}\cdot
\end{equation*}
\begin{equation}\label{SpC-GS2eq4}
\cdot\left(\pa\Delta_g S\pa_{2k}+(n-1)\kappa_0^2\pa S\pa_{2k}\right)^{\frac{1}{2}}\pa S\pa_{2k}^{k-\frac{1}{2}}+\pa S\pa_{2k}^k
\end{equation}
As we shall show in the Appendix C, it is impossible to give an upper bound to $\pa\Delta_g S\pa_{2k}$ in terms of $\lambda\pa S\pa_{2k}$. Nevertheless, we go beyond this difficulty using the following argument: we observe, as we are in the case where $S\in\mathcal A_1(\lambda)$ (resp. $S\in\mathcal A(\lambda)$), that $S$ and $\Delta_g S$ can be written as 
$$S=\sum_{i \mbox{ }s. t.\mbox{ }\lambda_i\leq\lambda} \alpha_i S_i\;,$$
$$\Delta_g S=\sum_{i\mbox{ }s.t.\mbox{ }\lambda_i\leq\lambda}\alpha_i\lambda_i S_i\;,$$

where $\{S_i\}$ is a $L^2$-orthonormal basis of eigenvectors for $\Delta_g$,
which implies that $\Delta_g S\in\mathcal A_1(\lambda)$ (resp. $\mathcal A(\lambda)$), and, if $A_p$ is the supremum of the ratio $\frac{\pa\phi\pa_p}{\pa\phi\pa_2}$ when $\phi$ runs in $\mathcal A_1(\lambda)\setminus\{0\}$ (resp. in $\mathcal A(\lambda)\setminus\{0\}$), we have:
\begin{equation}\label{SpC-GS2eq5}
\frac{\pa\Delta_g S\pa_{2k}}{\pa\Delta_g S\pa_2}\leq A_{2k}
\end{equation}
On the other hand the decomposition of $S$ in terms of the $S_i$'s and the Parseval identity give:
$$\pa\Delta_g S\pa_{2}^2\,=\sum_{i\mbox{ }s.t.\mbox{ }\lambda_i\leq\lambda}\lambda_i^2\alpha_i^2\leq\lambda^2\pa S\pa_2^2$$
This inequality and inequality (\ref{SpC-GS2eq5}) give
$$\frac{\pa\Delta_g S\pa_{2k}}{\pa S\pa_2}=\frac{\pa\Delta_g S\pa_{2k}}{\pa\Delta_g S\pa_2}\cdot\frac{\pa\Delta_g S\pa_{2}}{\pa S\pa_2}\leq A_{2k}\,\lambda$$
Now bearing this estimate in equation (\ref{SpC-GS2eq4}), we obtain:
$$A_{\frac{2kn}{n-2}}^k\leq \frac{[B(\alpha) D]\cdot k}{\sqrt{2k-1}}\cdot\left[ A_{2k}\,\lambda+(n-1)\,\kappa_0^2 \,A_{2k}\right]^{\frac{1}{2}}\cdot A_{2k}^{k-\frac{1}{2}}+A_{2k}^k$$
and thus
$$A_{\frac{2kn}{n-2}}\leq\left[1+\frac{[B(\alpha) D]\cdot k}{\sqrt{2k-1}}\left(\lambda+(n-1)\kappa_0^2\right)^{\frac{1}{2}}\right]^{\frac{1}{k}}\cdot A_{2k}$$
Now let us replace $k$ by $\beta^i$ (we recall that $\beta=\frac{n}{n-2}$) we see that
$$A_{2\beta^{m}}=\prod_{i=0}^{m-1}\frac{A_{2\beta^{i+1}}}{A_{2\beta^i}}\leq\prod_{i=0}^{m-1}\left[1+[B(\alpha)D]\cdot\frac{\beta^i}{\sqrt{2\beta^i-1}}(\lambda+(n-1)\kappa_0^2)^{\frac{1}{2}}\right]^{\frac{1}{\beta^i}}$$
Now, letting $m$ go to infinity, we obtain
$$A_\infty\leq\prod_{i=0}^\infty\left[1+[B(\alpha)D]\cdot\frac{\beta^i}{\sqrt{2\beta^i-1}}(\lambda+(n-1)\kappa_0^2)^{\frac{1}{2}}\right]^{\frac{1}{\beta^i}}$$
and we deduce that:
$$\frac{\pa S\pa_\infty}{\pa S\pa_2}\leq\xi\left([B(\alpha) D]\cdot(\lambda+(n-1)\kappa_0^2)^{\frac{1}{2}}\right).$$
When $S=f\in\mathcal A(\lambda)$ this becomes
$$\pa f\pa_\infty\leq\xi\left(B(\alpha) D\sqrt\lambda\right)\pa f\pa_2,$$
which proves the first inequality of (i);
when $S=df\in\mathcal A_1(\lambda)=d[\mathcal A(\lambda)]$ we obtain
$$\pa df\pa_\infty\leq\xi\left(B(\alpha) D\sqrt{\lambda +(n-1)\kappa^2}\right)\pa df\pa_2$$
which proves the first inequality of (ii). 
We conclude the proof by noticing that $\xi(x)\leq\exp\left(\frac{n}{2}\cdot\frac{x}{1+x}\right)(1+x)^{\frac{n}{2}}$, as proved in Appendix B. $\Box$\newline

\section{Spectral comparison between different manifolds in the presence of a map with bounded energy}

We find useful to introduce the following definitions:

\begin{defn}
Let $(Y,g)$, $(X,g_0)$ be two compact, connected Riemannian manifold. Let $F:(Y,g)\f (X,g_0)$ be a Lipschitz map. In every point $y\in Y$ where $F$ is differentiable (thus in almost every point of $Y$) we can define the pointwise energy of the map $F$ at $y$ as $$e_y(F)=\sum_1^n g_0(d_y F(e_i),d_y F(e_i))$$ where $\{e_i\}$ is any $g$-orthonormal basis of $T_yY$. Hence the global energy of the map $F$ is given by integration: $\mathrm{E}(F)=\int_Y e_y(F)\,dv_g(y)$.
\end{defn}

Throughout the paper we shall use the notion of \textit{absolute degree} of a continuous map $f$ between two $n$-dimensional compact manifolds ($\adeg(f)$) which is a homotopy invariant. Instead of defining here the absolute degree (the definition can be found, for example, in \cite{Epst}, \S 1) we shall give the notion of \textit{geometric degree} (more suitable for our purposes) and we remark that in \cite{Epst} Epstein proved that they are actually equal. Moreover, they coincide with the absolute value of the usual cohomological degree in case $f$ is a map between orientable manifolds.

\begin{defn}[Geometric degree]\label{SpC-gdeg}
Given a continuous map between two $n$-dimensional compact manifolds $f:Y\f X$, the geometric degree of $f$ is defined as $$\gdeg(f)=\inf\{G(h)\;|\;h:Y\f X\mbox{ is properly homotopic to } f\}$$
where $G(h)$ denotes the minimum number of connected components of\\ $h^{-1}(D)$, where $D$ varies among the top dimensional $n$-cells of $X$ such that $h: h^{-1}(D)\f D$ is a covering (if such a disk does not exist we say that $G(h)=\infty$).
\end{defn}

\subsection{A Bienaym\'e-\v{C}eby\v{s}\"ev inequality}

\begin{lemma}\label{SpC-BC}
Let $(Y,g)$ and $(X,g_0)$ be two connected, compact Riemannian manifolds of the same dimension, satisfying the following inequality between volumes:
$$\Vol_g(Y)\cdot(1-\eta)\leq\Vol_{g_0}(X),\quad\mbox{for some }\eta\in[0,1),$$
and let $F:(Y,g)\f(X,g_0)$ be a Lipschitz map with non-zero absolute degree, such that $|\jac(F)(y)|\leq(1+\eta)$, in every point $y$ where $F$ is differentiable. Let us define the set:
$Y_\eta^F=\{y\in Y|\, |\jac(F)(y)|\leq (1-\sqrt\eta)\}$. Then we have:
$$\frac{\Vol_{g}(Y_\eta^F)}{\Vol_g(Y)}\leq 2 \cdot\sqrt\eta.$$
\end{lemma}

\textbf{Proof.} Since the map $F$ is Lipschitz, it is differentiable almost everywhere, so the bounds that we gave in the statement are valid almost everywhere; we can apply the coarea formula (\cite{BZ}, Theorem 13.4.2) and we get:
$$\adeg(F)\cdot\Vol_{g_0}(X)\leq\int_X\#(F^{-1}(\{x\}))\,dv_{g_0}(x)=
$$
$$=\int_{Y}|\jac(F)(y)|\,dv_g(y)\leq(1-\sqrt\eta)\cdot\Vol_g(Y_\eta^F)+(1+\eta)\cdot\Vol_g(Y\setminus Y_\eta^F)$$
Since we are assuming $\adeg(F)\neq0$ and because of the inequality between the volumes of $(Y,g)$ and $(X,g_0)$ we obtain:
$$(1-\eta)\cdot\Vol_g(Y)\leq(1-\sqrt\eta)\cdot\Vol_g(Y_\eta^F)+(1+\eta)\cdot(\Vol_g(Y)-\Vol_g(Y_\eta^F))$$
thus we infer,
$$\frac{\Vol_g(Y_\eta^F)}{\Vol_g(Y)}\leq 2\cdot\sqrt\eta$$
which is the required inequality. $\Box$\newline

Let us consider a Lipschitz map $F:(Y,g)\f (X,g_0)$ between two compact, connected Riemannian manifolds, which satisfy the inequality  $$\Vol_g(Y)\cdot(1-\eta)\leq\Vol_{g_0}(X)\;.$$ If we assume that $\adeg(F)\neq 0$ and that the pointwise energy of the map $F$ satisfies the upper bound $e_y(F)\leq n\cdot(1+\eta)^{2/n}$, at almost every point $y$, then $F$ satisfies the assumptions of  Lemma \ref{SpC-BC}: in fact the geometric-arithmetic inequality $\prod_{i=1}^n\lambda_i^2\leq\left(\frac{1}{n}\lambda_i^2\right)^n$ implies that
\begin{equation*}|\det(d_yF)|\leq\left(\frac{1}{n}\tr((d_yF)^t\circ(d_yF))\right)^{\frac{n}{2}}= 
\end{equation*}
\begin{equation}\label{SpC-detEST}
=\left[\frac{1}{n}\, e_y(F)\right]^{\frac{n}{2}}\leq1+\eta
\end{equation}
Under these new assumptions we obtain the following

\begin{lemma}\label{SpC-qi} Let $0<\eta\leq\frac{1}{4}$. 
Let $(X,g_0)$, $(Y,g)$ be two compact, connected Riemannian manifolds which satisfy the inequality $\Vol_g(Y)(1-\eta)\leq\Vol_{g_0}(X)$. In any point $y\in Y$ such that $|\jac(F)(y)|\geq(1-\sqrt\eta)$ and $e_y(F)\leq n(1+\eta)^{\frac{2}{n}},$ $d_y F$ is a quasi-isometry; more precisely we have $\forall u\in T_yY$
$$(1-5(n-1)\eta^{\frac{1}{4}})\pa u\pa_g^2\,\leq\,\pa d_yF(u)\pa_{g_0}^2\,\leq\,(1+5(n-1)\eta^{\frac{1}{4}})\pa u\pa_g^2$$
\end{lemma}

\textbf{Proof.} Let us consider the bilinear symmetric form given by
$$(u,v)\f g_0(d_y F(u),d_yF(v))=g((d_yF)^t\circ(d_yF)(u), v)$$
defined on $T_yY\times T_yY$. We denote by $A$ the matrix associated to the endomorphism $(d_yF)^t\circ(d_yF)$ in a $g$-orthonormal basis of $(T_yY, g_y)$; this matrix being symmetric and non-negative with determinant greater or equal to $(1-\sqrt\eta)^2$ by assumption and with trace equal to $e_y(F)$ (hence, by assumption, less or equal to $n\,(1+\eta)^{\frac{2}{n}}$), we can apply Proposition A.1 in Appendix A, which gives, $\forall u\in T_yY$,
$$|\pa d_yF(u)\pa_{g_0}^2-\pa u\pa_{g}^2|=|g((A-\Id)u,u)|\leq$$
$$\leq 2(n-1)\,\eta^{\frac{1}{4}}\left(1+\frac{n+10}{2n}\sqrt\eta\right)^{\frac{1}{2}}\pa u\pa_g^2\,\leq\, 5(n-1)\eta^{\frac{1}{4}}\pa u\pa_g^2\quad \Box$$

\subsection{A general comparison Lemma}

Let us start  with some definitions:

\begin{defn}\label{SpC-rayleigh} 
Let $f$ be a function in $H_1^2(M,g)\setminus\{0\}$ where $(M,g)$ is a fixed Riemannian manifold; the Rayleigh quotient of $f$ is defined as the positive real number:
$$\mathcal R_g(f)=\frac{\int_M|df|^2\,dv_g}{\int_M |f|^2\,dv_g}.$$
\end{defn}

\begin{defn} We will denote by $\lambda_i(X,g_0)$ and $\lambda_i(Y,g)$ the eigenvalues of the Laplace-Beltrami operators of $(X,g_0)$ and $(Y,g)$ respectively, indexed in increasing order and counted with their multiplicity from zero to infinity. 
\end{defn}

We denote by $\mathcal A_X(\lambda)$ the direct sum of the eigenspaces of the Laplace-Beltrami operator of $(X,g_0)$ corresponding to the eigenvalues which are less or equal to $\lambda$. We remark that the eigenvalues of  the Laplacian $\Delta_g$ of $(Y,g)$, are also the eigenvalues of the quadratic form $u\f\int_Y |du|^2\,dv_g$ with respect to the $L^2$-scalar product. A classical consequence of the Minimax Principle is the
\begin{lemma}\label{SpC-comparison}
If there exists a linear map $\phi:\mathcal A_X(\lambda)\f H_1^2(Y)$ such that $\forall u\in\mathcal A_X(\lambda)$
\begin{itemize}
\item[(i)] $\pa\phi(u)\pa_2^2\geq(1-\delta)\pa u\pa_2^2$ where $\delta\in[0,1)$;
\item[(ii)] $\frac{1}{\Vol_g(Y)}\int_{Y}|d(\phi(u))|^2\,dv_g\,\leq\,(1+\varepsilon)\frac{1}{\Vol_{g_0}(X)}\int_X|du|^2\,dv_{g_0}$;
\end{itemize}
then, for every $i\in\N$ such that $\lambda_i(X,g_0)\leq\lambda$, we have $$\lambda_i(Y,g)\leq\left(\frac{1+\varepsilon}{1-\delta}\right)\lambda_i(X,g_0).$$
\end{lemma}

\subsection{Proof of the Proposition \ref{SpC-tech1}}
We start with the following lemma:

\begin{lemma}\label{SpC-tech2}
Under the assumptions of Proposition \ref{SpC-tech1}, for any function $f:X\f\R$ such that $f\in\mathcal A_X(\lambda)$, if $\eta\leq\frac{1}{9}$ we have:
$$\mathrm{(i)}\mbox{ }\left(\frac{1-\eta}{1+\eta}\right)\cdot\frac{1}{\Vol_{g_0}(X)}\int_X f^2\,dv_{g_0}\leq\frac{1}{\Vol_g(Y)}\int_Y(f\circ F)^2\,dv_g\leq$$
$$\leq\left[1+3\sqrt\eta\left(1+\xi^2(B(\alpha)D\sqrt\lambda)\right)\right]\cdot\frac{1}{\Vol_{g_0}(X)}\int_X f^2\,dv_{g_0}$$
$$\mathrm{(ii)}\mbox{  }\frac{1}{\Vol_g(Y)}\int_Y|d(f\circ F)|^2\,dv_g\leq(1+5(n-1)\eta^{\frac{1}{4}})\cdot$$
$$\cdot\left(1+\left[2+(2n+1)\xi^2(B(\alpha)\sqrt{\lambda D^2+(n-1)\alpha^2})\right]\eta^{\frac{1}{2}}\right)\cdot\frac{1}{\Vol_{g_0}(X)}\int_X|df|^2\,dv_{g_0}$$
\end{lemma}

\textbf{Proof.} We start with the proof of the first inequality of the property (i). As discussed in section 3.1 (see equation (\ref{SpC-detEST})), the bound on the pointwise energy of $F$ implies that $|\jac(F)|$ is bounded above by $(1+\eta)$. Hence, using the assumptions and the coarea formula (\cite{BZ}, Theorem 13.4.2), we find:
$$\frac{1+\eta}{\Vol_g(Y)}\cdot\int_Y |(f\circ F)(y)|^2\,dv_g(y)\geq$$
$$\geq\frac{1}{\Vol_{g}(Y)}\cdot\int_Y|(f\circ F)(y)|^2\;|\jac(F)(y)|\,dv_g(y)\geq$$
$$\geq\frac{1-\eta}{\Vol_{g_0}(X)}\cdot\int_X|f(x)|^2\,\#(F^{-1}(\{x\}))\,dv_{g_0}(x)\geq\frac{1-\eta}{\Vol_{g_0}(X)}\cdot\int_X|f(x)|^2\,dv_{g_0}(x)$$
which proves the first inequality of (i). We remark that the second inequality of  (i) is not necessary in order to prove the Proposition \ref{SpC-tech1}. However we shall provide a proof of this inequality for the sake of completeness. We notice that we are under the assumptions of the Lemma \ref{SpC-BC} and thus $\frac{\Vol_g(Y_\eta^F)}{\Vol_g(Y)}\leq2\,\sqrt\eta$ where $Y_\eta^F=\{y\in Y\;\;|\;\;|\jac(F)|(y)<1-\sqrt\eta\}$. 
As the absolute degree is not trivial, $F$ is surjective and thus $\#(F^{-1}(\{x\}))\geq1$ for every $x\in X$.
Thus we can get a first estimate for the $L^2$-norm of $F^*(f)=f\circ F$ in terms of the $L^\infty$-norm and the $L^2$-norm of $f$; in fact using the coarea formula (\cite{BZ}, Theorem 13.4.2) we find  on $Y\setminus Y_\eta^F$
\begin{equation*}
(1-\sqrt\eta)\cdot\int_{Y\setminus Y_\eta^F}|f\circ F|^2\,dv_g\leq \int_Y|f\circ F|^2|\jac(F)(y)|\,dv_g=
\end{equation*}
\begin{equation*}
=\int_X\#(F^{-1}(\{x\}))|f(x)|^2\,dv_{g_0}(x)\leq
\end{equation*}
\begin{equation}\label{SpC-**}
\leq\Vol_{g_0}(X)\pa f\pa_2^2+\pa f\pa_\infty^2\int_X(\#(F^{-1}(\{x\}))-1)\,dv_{g_0}
\end{equation}
On the other hand we know that $\Vol_g(Y)(1-\eta)<\Vol_{g_0}(X)$ and from the coarea formula and the upper bound on the Jacobian of $F$ we deduce:
\begin{equation*}
0\leq\int_X(\#(F^{-1}(\{x\}))-1)\,dv_{g_0}\leq
\end{equation*}
\begin{equation}\label{SpC-***}
\leq\int_Y|\jac(F)(y)|\,dv_g(y)-(1-\eta)\Vol_g(Y)\leq2\eta\cdot\Vol_g(Y)
\end{equation}
using the fact that $|\jac(F)|\leq1+\eta$ \textit{a.e.} and that $F$ is surjective and applying the coarea formula, we get 
\begin{equation}\label{SpC-numeroter}
\Vol_{g_0}(X)\leq\int_Y|\jac(F)(y)|\,dv_g(y)\leq(1+\eta)\Vol_g(Y)
\end{equation}
hence we obtain, from (\ref{SpC-**}), (\ref{SpC-***}) and (\ref{SpC-numeroter}):
$$\int_{Y\setminus Y_\eta^F}|f\circ F|^2\,dv_{g}\leq\left(\frac{1+\eta}{1-\sqrt\eta}\pa f\pa_2^2+\frac{2\eta}{1-\sqrt\eta}\pa f\pa_\infty^2\right)\cdot\Vol_g(Y)$$
whereas on $Y_\eta^F$, using Lemma \ref{SpC-BC}, we infer
$$\int_{Y_\eta^F}|f\circ F|^2\,dv_g\leq\pa f\pa_\infty^2\cdot\Vol_g(Y_\eta^F)\leq2\sqrt\eta\cdot\pa f\pa_\infty^2\cdot\Vol_g(Y)$$
we sum these two inequalities, and we divide both sides by $\Vol_g(Y)$:
$$\pa f\circ F\pa_2^2\leq\left(\frac{1+\eta}{1-\sqrt\eta}\right)\cdot\left(\pa f\pa_2^2+\frac{2\sqrt\eta}{1+\eta}\pa f\pa_\infty^2\right)$$
Now we can use Proposition \ref{SpC-GSobolev2} (i), which tells us that, for every $f\in\mathcal A_X(\lambda)$,
$$\frac{\pa f\pa_\infty^2}{\pa f\pa_2^2}\leq\xi^2(B(\alpha) D\sqrt\lambda)$$
where $\alpha=\kappa D$. Thus we obtain the estimate:
$$\pa f\circ F\pa_{L^2(Y)}^2\leq\left(\frac{1+\eta}{1-\sqrt\eta}\right)\left[1+\frac{2\sqrt\eta}{1+\eta}\cdot\xi^2(B(\alpha) D\sqrt\lambda)\right]\cdot\pa f\pa_{L^2(X)}^2$$
which can be simplified, thanks to the assumption $\eta\leq\frac{1}{9}$ in
$$\pa f\circ F\pa_{L^2(Y)}^2\leq\left[1+3\sqrt\eta\left(1+\xi^2(B(\alpha)D\sqrt\lambda)\right)\right]\cdot\pa f\pa_{L^2(X)}^2$$
This ends the proof of inequalities (i) of the Lemma \ref{SpC-tech2}.\\
Now we shall prove the inequality (ii). For every $y\in Y$ we have
$$|d_y(f\circ F)|^2=\sup_{u\in T_yY\setminus\{ 0_y\}}\left(\frac{|d_y(f\circ F)(u)|^2}{| u|^2}\right)\leq$$
$$\leq\sup_{u\in T_yY\setminus\{0_y\}}\left(\frac{|d_{F(y)}f(d_yF(u))|^2} {|d_yF(u)|^2}\cdot\frac{|d_yF(u)|^2}{|u|^2}\right)$$
and thus
$$|d_y(f\circ F)|^2\leq|d_{F(y)}f|^2\cdot\pa\mid d_yF\mid\pa^2$$
where $\pa\mid d_y F\mid\pa$ denotes the operator norm of $d_yF$.
By  Lemma \ref{SpC-qi} we get for every $y$ in $ Y\setminus Y_\eta^F$
\begin{equation}\label{SpC-LA1}
|d_y(f\circ F)|^2\leq(1+5(n-1)\eta^{\frac{1}{4}})\cdot|d_{F(y)}f|^2
\end{equation}
whereas for every $y\in Y_\eta^F$ we have $\pa\mid d_yF\mid\pa^2\leq e_y(F)\leq n(1+\eta)^{\frac{2}{n}}$, which gives $\forall y\in Y_\eta^F$
\begin{equation}\label{SpC-LA2}
|d_y(f\circ F)|^2\leq n(1+\eta)^{\frac{2}{n}}\cdot|d_{F(y)}f|^2\leq n(1+\eta)^{\frac{2}{n}}\pa df\pa_\infty^2
\end{equation}
From equation (\ref{SpC-LA2}) we deduce that:
\begin{equation*}
\frac{1}{\Vol_g(Y)}\int_{Y_\eta^F}|d(f\circ F)|^2\,dv_g\leq
\end{equation*}
\begin{equation}\label{SpC-LA3}
\leq\frac{\Vol_g(Y_\eta^F)}{\Vol_g(Y)}\,n\,(1+\eta)^{\frac{2}{n}}\pa df\pa_\infty^2\leq
2n\,\sqrt\eta\, (1+\eta)^{\frac{2}{n}}\pa df\pa_\infty^2\;,
\end{equation}
where the last inequality deduces from the Lemma \ref{SpC-BC}.
On the other hand, from equation (\ref{SpC-LA1}) and from the definition of $Y_\eta^F$, we deduce that
$$\frac{1}{\Vol_g(Y)}\int_{Y\setminus Y_\eta^F}|d(f\circ F)|^2\,dv_g\leq$$
$$\leq\frac{(1+5(n-1)\eta^{\frac{1}{4}})}{\Vol_g(Y)(1-\sqrt\eta)}\int_{Y\setminus Y_\eta^F} |d_{F(y)}f|^2|\jac(F)(y)|\,dv_g(y)\leq$$
$$\leq\frac{(1+5(n-1)\eta^{\frac{1}{4}})}{\Vol_g(Y)(1-\sqrt\eta)}\int_X|d_xf|^2(\#(F^{-1}(\{ x\})))\,dv_{g_0}(x)\leq$$
$$\leq\frac{(1+5(n-1)\eta^{\frac{1}{4}})}{\Vol_g(Y)(1-\sqrt\eta)}\left(\int_X |df|^2\,dv_{g_0}+\pa df\pa_\infty^2\int_X(\# F^{-1}(\{x\})-1)\,dv_{g_0}(x)\right),$$
where the second inequality follows from the coarea formula (\cite{BZ}, Theorem 13.4.2). From this inequality and from inequalities (\ref{SpC-numeroter}) and (\ref{SpC-***}) we deduce
$$\frac{1}{\Vol_g(Y)}\int_{Y\setminus Y_\eta^F}|d(f\circ F)|^2\,dv_g\leq$$
$$\leq\frac{(1+5(n-1)\eta^{\frac{1}{4}})}{1-\sqrt\eta}\cdot\left[\frac{\Vol_{g_0}(X)}{\Vol_g(Y)}\pa df\pa_2^2+2\eta\pa df\pa_\infty^2\right]$$
Now we sum the last inequality with equation (\ref{SpC-LA3}); 
and we obtain, using the fact that $\eta\leq\frac{1}{9}$:
$$\pa d(f\circ F)\pa_2^2\leq\frac{(1+5(n-1)\eta^{\frac{1}{4}})}{(1-\sqrt\eta)}\left[(1+\eta)\pa df\pa_2^2+2\eta\pa df\pa_\infty^2\right]+$$
$$+2\,n\sqrt\eta\,(1+\eta)^{\frac{2}{n}}\pa df\pa_\infty^2$$
$$\leq(1+5(n-1)\eta^{\frac{1}{4}})\left[(1+2\sqrt\eta)\pa df\pa_2^2+(2n+1)\sqrt\eta \pa df\pa_\infty^2\right]$$
To conclude it is sufficient to apply Proposition \ref{SpC-GSobolev2} (ii) which gives
$$\pa df\pa_\infty^2\leq\xi^2(B(\alpha)\sqrt{\lambda D^2+(n-1)\alpha^2})\pa df\pa_2^2$$
and thus achieves the proof of inequality (ii) of the Lemma \ref{SpC-tech2}.\;$\Box$\\

\textbf{End of the Proof of the Proposition \ref{SpC-tech1}.} We just apply the Lemma \ref{SpC-comparison} to the linear map $F^*:f\f f\circ F$. This map is linear and sends $\mathcal A_X(\lambda)$ onto a subspace of $H_1^2(Y,g)$: in fact, as $f$ is $C^\infty$, $f\circ F$ is continuous and Lipschitz, thus $f\circ F$ and $|d(f\circ F)|$ are bounded and have finite $L^2$-norms, this proves that $f\circ F\in H_1^2(Y,g)$ and that $F^*\mathcal [A_X(\lambda)]$ is included in $H_1^2(Y,g)$. By the Lemma \ref{SpC-tech2} the assumptions of Lemma \ref{SpC-comparison} are verified for every $f\in\mathcal A_X(\lambda)$ and, applying the Lemma \ref{SpC-comparison} , we obtain
$$\lambda_i(Y,g)\leq\left(\frac{(1+5(n-1)\eta^{\frac{1}{4}})(1+\eta)}{1-\eta}\right)\cdot$$
$$\cdot\left(1+\left[2+(2n+1)\xi^2\left(B(\alpha)\sqrt{\lambda D^2+(n-1)\alpha^2}\right)\right]\eta^{\frac{1}{2}}\right)\lambda_i(X,g_0).$$
Using the estimate of $\xi$ computed in Appendix B we see that:
$$\lambda_i(Y,g)\leq(1+7(n-1)\eta^{\frac{1}{4}})\cdot$$
$$\cdot\left(1+\left[2+(2n+1)e^{n}\left(1+B(\alpha)\sqrt{\lambda D^2+(n-1)\alpha^2}\right)^n\right]\eta^{\frac{1}{2}}\right)\cdot\lambda_i(X,g_0)$$
We conclude by taking $\lambda=\lambda_i(X,g_0)$ in the last inequality. $\Box$

\section{Spectral comparison between manifolds in terms of their Gromov-Hausdorff distance}

The main purpose of this section is to present the link between the spectra comparison theorem which we proved in the previous section and the barycenter method by Besson, Courtois and Gallot (see \cite{BCG1},\cite{BCG2}). More precisely we will use a recent developement of this technique by L. Sabatini, \cite{Sab}. The main feature of this last version of the barycenter method is that, on one hand, no assumption is made on the sign of the sectional curvature of the \lq\lq known" manifold $(X,g_0)$ (only its boundedness is required), on the other hand, no condition is assumed on the geometry of the \lq\lq unknown" manifold $(Y,g)$, except for the fact that the Gromov-Hausdorff distance between $(Y,g)$ and $(X,g_0)$ is supposed to be smaller than some universal constant $\varepsilon_0$, which is precised in the Theorem \ref{SpC-barycenter}. This technique, combined with Proposition \ref{SpC-tech1}, will provide a spectra comparison theorem between manifolds satisfying weak assumptions. 

\subsection{Proof of Theorem \ref{SpC-MT2}} In his PhD thesis L. Sabatini proved the following theorem:

\begin{thm}[L. Sabatini, \cite{Sab}]\label{SpC-barycenter}
Let $(X,g_0)$ be a compact Riemannian manifold of dimension $n$, whose sectional curvature $\sigma$ satisfies the bound $|\sigma|\leq K^2$, for a suitable $K>0$. Let $i_0$ denote the injectivity radius of $(X,g_0).$ Let $(Y,g)$ be another compact Riemannian manifold such that there exists a measurable Gromov-Hausdorff $\varepsilon$-approximation $f:(Y,g)\f(X,g_0)$ with $\varepsilon<\varepsilon_0=\varepsilon_0(n,i_0,K)$
(an explicit value for $\varepsilon_0$ can be found in \cite{Sab}, Theorem 4.4.1) , then there exists a $C^1$-map $F:(Y,g)\f(X,g_0)$ with the following properties:
\begin{enumerate}
\item For any $y\in Y$ one has:
$$e_y(F)\leq n\left(1+20(n+1)(K \varepsilon)^{\frac{1}{4}}\right)\;,$$
$$|\jac(F)|\leq\left(1+20(n+1)(K\varepsilon)^{\frac{1}{4}}\right)^{\frac{n}{2}}\;.$$
\item If moreover the $\varepsilon$-Hausdorff approximation is continuous the map $F$ is homotopic to $f$.
\end{enumerate}
\end{thm}
L. Sabatini used this result to provide a sharp lower bound to the ratio $\frac{\Vol_g(Y)}{\Vol_{g_0}(X)}$ in terms of $\varepsilon$. Let us see how to use this theorem in order to end the proof of Theorem \ref{SpC-MT2}.\\

\textbf{End of the proof of Theorem \ref{SpC-MT2}.}  Any pair of Riemannian manifolds which satisfies the assumptions of Theorem \ref{SpC-MT2} also satisfies the assumption of Theorem \ref{SpC-barycenter}. Applying Theorem \ref{SpC-barycenter} (2), we obtain the existence of a $C^1$ map $F:(Y,g)\f(X,g_0)$ homotopic to $f$ and thus of non zero absolute degree which (by Theorem \ref{SpC-barycenter} (1)) satisfies, at every point $y\in Y$,
$$e_y(F)\leq n\,(1+\eta(\varepsilon))^{\frac{2}{n}}\;,$$
where $\eta(\varepsilon)$ is defined by
$$\eta(\varepsilon)=\left[1+20\,(n+1)\,(\kappa\,\varepsilon)^{\frac{1}{4}}\right]^{\frac{n}{2}}-1\;.$$
The assumption $\varepsilon<\frac{1}{\kappa}\left(\frac{\left(\frac{10}{9}\right)^{\frac{2}{n}}-1}{20\,(n+1)}\right)^4$ immediately implies that $\eta(\varepsilon)<\frac{1}{9}$. Finally, the assumption:
$$\frac{\Vol_{g_0}(X)}{\Vol_g(Y)}\geq 1-10\,n\,(n+1)\cdot(\kappa\,\varepsilon)^{\frac{1}{4}}$$
implies that $\frac{\Vol_{g_0}(X)}{\Vol_g(Y)}>1-\eta(\varepsilon)$ because $(1+x)^{\frac{n}{2}}-1\,\geq\,\frac{n}{2}\, x$, $\forall x\in\R^+$. We may thus apply the Proposition \ref{SpC-tech1} which proves that 
$$\frac{\lambda_i(Y,g)}{\lambda_i(X,g_0)}\leq\left[1+7\,(n-1)\,\eta(\varepsilon)^{\frac{1}{4}}\right]\cdot\left[1+C(n,\,D^2\,\lambda_i(X,g_0),\,\kappa\, D)\,\eta(\varepsilon)^{\frac{1}{2}}\right]$$
we conclude the proof when noticing that, $\forall x\in\R^+$ we have $(1+x)^{\frac{n}{2}}-1\,\leq\,\frac{n}{2}\,(1+x)^{\frac{n}{2}}\,x$, and thus
$$\eta(\varepsilon) \leq 10\,n\,(n+1)\,(1+\eta(\varepsilon))\,(\kappa\,\varepsilon)^{\frac{1}{4}}\leq\frac{100}{9}\,n\,(n+1)\,(\kappa\,\varepsilon)^{\frac{1}{4}}$$
which leads to
$$7\,(n-1)\,\eta(\varepsilon)^{\frac{1}{4}}\,<\,14\,(n-1)\,\sqrt{n+1}\,(\kappa\varepsilon)^{\frac{1}{16}}=C_1(n)\,(\kappa\varepsilon)^{\frac{1}{16}}\;,$$
$$C(n,\,D^2\,\lambda_i(X,g_0),\,\kappa\, D)\,\eta(\varepsilon)^{\frac{1}{2}}\,\leq\,\frac{10}{3}\,(n+1)\,(\kappa\,\varepsilon)^{\frac{1}{8}}\cdot$$
$$\cdot\left(2+2(n+1)\,e^n\,\left[1+B(\kappa\,D)\cdot D\cdot\sqrt{\lambda_i(X,g_0)+(n-1)\kappa^2}\right]^n\right)\leq$$
$$\leq C_2(n,\,\kappa\, D,\,D^2\,\lambda_i(X,g_0))\cdot(\kappa\,\varepsilon)^{\frac{1}{8}}$$
this concludes the proof of the Theorem \ref{SpC-MT2}. $\Box$\\

\subsection{Examples} This subsection is devoted to the construction of examples and counterexamples regarding Theorem \ref{SpC-MT2}.

\begin{exmp}\label{SpC-4.2}
Consider any closed Riemannian manifold $(X^n,g_0)$. Let us fix $(X,g_0)$ and call $D,\,i_0,\,\kappa$ its diameter, its injectivity radius and the maximum of $|\sigma_0|^{1/2}$ (where $\sigma_0$ denotes the sectional curvature of $(X,g_0)$). Starting from $(X,g_0)$ we shall construct a family $\{(Y_\varepsilon,g_\varepsilon)\}_{\varepsilon>0}$ of Riemannian manifolds with the following properties:
\begin{enumerate}
\item There exists a continuous map $f:(Y_\varepsilon,g_\varepsilon)\f(X,g_0)$ of nonzero absolute degree which is a Gromov-Hausdorff $\varepsilon$-approximation for $\varepsilon< \varepsilon_1(n,i_0,\kappa)$.
\item The assumption $[1-10\,n\,(n+1)(\kappa\,\varepsilon)^{\frac{1}{4}}]\,\Vol{g_\varepsilon}(Y_\varepsilon)<\Vol_{g_0}(X)$ is satisfied
(actually we shall construct a sequence of manifolds satisfying the stronger condition $\Vol_{g_\varepsilon}(Y_\varepsilon)<\Vol_{g_0}(X)$).
\end{enumerate}
Thus, for every $\varepsilon>0$, the pair of Riemannian manifolds $\{(X,g_0)\,;\,(Y_\varepsilon,g_\varepsilon)\}$ satisfies all the assumptions of the Theorem \ref{SpC-MT2}.\\\\
Let $x_0\in X$. Let $\rho(\varepsilon)=\frac{\varepsilon}{4}$, we excide the geodesic ball $B(x_0,\rho(\varepsilon))$ (we remark that $\varepsilon<\varepsilon_1(n,\,i_0,\,\kappa)<\inj(X,g_0)$). Consider any compact, $n$-dimensional Riemannian manifold $(Z_\varepsilon,h)$. We rescale the metric $h$ by a scale factor $$\alpha(\varepsilon)^2=\min\left\{\frac{\varepsilon^2}{16\,[\diam(Z,h)]^2}\;;\,\left(\frac{\Vol_{g_0}(B(x_0,\rho(\varepsilon)))}{2\,\Vol_h(Z_\varepsilon)}\right)^{\frac{2}{n}}\right\}$$
we shall refer to the rescaled metric as $h_\varepsilon$. Let us excide a  ball $B(z_0,r_\varepsilon)$ of radius $r_\varepsilon=\frac{\inj(Z_\varepsilon,h_\varepsilon)}{2}$ (which is strictly less than $\varepsilon/4$, by construction). We glue $(Z_\varepsilon,h_\varepsilon)\setminus B(z_0,r_\varepsilon)$ and $(X,g_0)\setminus B(x_0,\rho(\varepsilon))$ along a tube $S^{n-1}\times I$ endowed with a metric $k_\varepsilon$ which is given by
$$(k_\varepsilon)|_{(x,r)}=\left(
\begin{array}{cc}
r\cdot(h_\varepsilon)|_x+(1-r)\cdot g|_x &  0 \\
0 & l(\varepsilon)^2\\
\end{array}\right)$$
where
$$l(\varepsilon)=\min\left\{\frac{\varepsilon}{4}\;;\;\frac{\Vol_{g_0}(B(x_0,\rho(\varepsilon)))}{2\cdot\int_0^1\int_{S^{n-1}}\left[\det(r\cdot(h_\varepsilon)|_x+(1-r)\cdot g|_x)\right]^{\frac{1}{2}}\,dx\,dr}\right\}$$
Now we consider the resulting metric $g_\varepsilon$ over $Y_\varepsilon=Z_\varepsilon\#(S^{n-1}\times I)\# X$. We remark that this metric is only $C^\infty$ piecewise, being only continuous at the gluing spheres; however, by the choices made during the construction, it is clear that these spheres possess a tubular neighbourhood which is diffeomorphic to $S^{n-1}\times I$. Hence, using mollifiers on very small tubular neighbourhood of the gluing spheres  we can smooth the metric $g_\varepsilon$, without significant changes for the volume and the diameter. We shall call $g_\varepsilon$ the new (smooth) metric.
In particular we can arrange things in order to have that $(Y_\varepsilon,g_\varepsilon)$ satisfies the volume assumption of Theorem \ref{SpC-MT2}.\\\\
Now we define a map $f:Y_\varepsilon\f X$ by sending $Y_\varepsilon\setminus [Z_\varepsilon\# (S^{n-1}\times I)]$ identically on $X\setminus B(x_0,\rho(\varepsilon))$, the tube $S^{n-1}\times [0,1)$ on $B(x_0,\rho(\varepsilon))\setminus\{x_0\}$ (sending $(x,r)$ in $\exp_{x_0}((1-r)\cdot\rho(\varepsilon)\,\exp_{x_0}^{-1}(x))$) and $Z_\varepsilon\setminus B(z_0,r_\varepsilon)$ on $x_0$. This map has non zero absolute degree and is a Gromov-Hausdorff $\varepsilon$-approximation. To see that it is a $\varepsilon$-approximation observe that the map is surjective, hence it is sufficient to verify that for any $y,\,y'\in Y_\varepsilon$
$$|d_{g_\varepsilon}(y,y')-d_{g_0}(f(y),f(y'))|<\varepsilon$$
Let us first show the inequality $d_{g_\varepsilon}(y,y')<d_{g_0}(f(y),f(y'))+\varepsilon$; if $y,\,y'\in Z_\varepsilon\#(S^{n-1}\times I)$ this follows directly from the fact that $$\diam(Z_\varepsilon\#(S^{n-1}\times I),g_\varepsilon|_{Z_\varepsilon\#(S^{n-1}\times I)})<\frac{3\,\varepsilon}{4}\;.$$
 Assume now $y\in Z_\varepsilon\#(S^{n-1}\times I)$ and $y'\in Y_\varepsilon\setminus[Z_\varepsilon\#(S^{n-1}\times I)]$; take a length minimizing path from $f(y')$ to $f(y)$ and observe that it must meet 
$\de B(x_0,\rho(\varepsilon))$ at some point $x_1$; let us call $\gamma$ the path from $y'$ to $y_1$ (the point in $Y$ corresponding to the point where the  geodesic from $f(y')$ to $f(y)$ meets $\de B(x_0,\rho(\varepsilon))$ for the first time) and compose this path with a minimizing geodesic $\delta$ from $y_1$ to $y$, then:
$$d_{g_\varepsilon}(y,y')\leq\mathrm{length}(\gamma)+\mathrm{length}(\delta)< d_{g_0}(f(y), f(y'))+\frac{3\,\varepsilon}{4}$$
Finally, if both $y$ and $y'$ lay in $Y_\varepsilon\setminus[Z_\varepsilon\#(S^{n-1}\times I)]$, it is clear that 
$$d_{g_\varepsilon}(y,y')\leq d_{g_0}(f(y),f(y'))+\diam(Z_\varepsilon\#(S^{n-1}\times I),g_\varepsilon|_{Z\#(S^{n-1}\times I)})<$$
$$<d_{g_0}(f(y),f(y'))+\frac{3\,\varepsilon}{4}$$
Now we prove the inequality $d_{g_0}(f(y),f(y'))<d_{g_\varepsilon}(y,y')+\varepsilon$; first take $y,\,y'\in Z_\varepsilon\#(S^{n-1}\times I)$, then $f(y)$, $f(y')$ are both in $B(x_0,\rho(\varepsilon))$, hence their distance is less than $2\rho(\varepsilon)<\frac{\varepsilon}{2}$ and $d_{g_0}(f(y),f(y'))< d_{g_\varepsilon}(y,y')+\varepsilon$. Now assume that $y\in Z_\varepsilon\#(S^{n-1}\times I)$ and $y'\in Y_\varepsilon\setminus[ Z_\varepsilon\#(S^{n-1}\times I)]$; take a minimizing geodesic from $y'$ to $y$; this geodesic meets $\de B(x_0,\rho(\varepsilon))$ at a first point $y_1$. Call $\gamma$ the geodesic segment from $y'$ to $y_1$. Since $f(y)\in B(x_0,\rho(\varepsilon))$ which has diameter less than $\frac{\varepsilon}{2}$ we have that:
$$d_{g_0}(f(y),f(y'))<\mathrm{length}(\gamma)+\frac{\varepsilon}{2}\leq d_{g_\varepsilon}(y,y')+\frac{\varepsilon}{2}$$
Finally, if both $y$ and $y'$ lay outside $Z_\varepsilon\#(S^{n-1}\times I)$,  consider a minimizing geodesic from $y$ to $y'$; if the geodesic does not cross $\de B(x_0,\rho(\varepsilon))$ then it is a geodesic also for the metric $g_0$ on $X$. Otherwise we take the two geodesic segments of the minimizing geodesic joining $y$ and $y'$ respectively with $B(x_0,\rho(\varepsilon))$ and we join them by two geodesic rays centered at $x_0$. We call $\gamma$ this path; then, by construction we have:
$$d_{g_0}(f(y),f(y'))\leq d_{g_\varepsilon}(y,y')+\frac{\varepsilon}{2}\;.$$
This proves that $f$ is a Gromov-Hausdorff $\varepsilon$-approximation.\\\\
Let us finally point out that in the previous construction we did not make any topological assumption on $Z_\varepsilon$, thus this example shows that, for any $(X,g_0)$ compact Riemannian manifold, there are infinitely many pairwise non homotopic Riemannian manifolds $(Y,g)$ which satisfy the assumptions of Theorem \ref{SpC-MT2}.

\end{exmp}

\begin{exmp}[Mushrooms]\label{SpC-4.3}
We shall show the necessity of the assumption that we made on the volumes of the manifolds in Theorem \ref{SpC-MT2}, by adding 'mushrooms' to a fixed manifold. First we show how to construct a mushroom. We excise from $(X,g_0)$ a geodesic ball $B(x_0,\frac{\varepsilon}{2})$ of radius $\frac{\varepsilon}{2}$, where $\varepsilon<\frac{\varepsilon_1(n,i_0,\kappa)}{2}$. Let $V_{\varepsilon}=\Vol_{g_0}(B(x_0,\frac{\varepsilon}{2}))$. We shall glue to $(X,g_0)$ a standard sphere of radius $f(\varepsilon)$ such that:
$$f(\varepsilon)=\min\left\{\left(\frac{V_\varepsilon}{2}\right)^{\frac{1}{n}}\frac{1}{\omega_n};\;\frac{\varepsilon}{2\pi}\right\}\;.$$
Let $\delta<<f(\varepsilon)$ (afterwards we shall consider $\delta\f0$) and excise from $S^{n}(f(\varepsilon))$ a geodesic ball $B(z_0,\arcsin(\delta))$ from $S^n(f(\varepsilon))$.
We glue along $\de B(x_0,\frac{\varepsilon}{2})$ and $\de B(z_0,\arcsin(\delta))$ respectively a tube $S^{n-1}\times I$ endowed with a metric $g^{\delta,\varepsilon}$ defined as follows:
$$g^{\delta,\varepsilon}=[(1-\lambda^{\varepsilon}(r))\cdot g_0|_{B(x_0,(1-r)\frac{\varepsilon}{2})}+\lambda^{\varepsilon}(r)\cdot(f(\varepsilon))^2\cdot\can|_{\de B(z_0,r\arcsin(\delta)))}]\oplus\delta^2\cdot dr^2$$
where '$\can$' is the standard metric of $S^{n}$ and
where $\lambda^{\varepsilon}$ is an increasing $C^\infty$ function defined on $I=[0,1]$ satisfying:
$$\lambda^{\varepsilon}|_{[0,\frac{1}{3})}=0,\;\lambda^{\varepsilon}|_{(\frac{2}{3},1]}=1\,;$$
we choose $\delta$ sufficiently small that we have:
$$\delta<\frac{\varepsilon}{4}\;;\quad\Vol_{g^{\delta,\varepsilon}}(S^{n-1}\times I)<\frac{V_{\varepsilon}}{2}\;.$$
In particular the resulting Riemannian manifold is diffeomorphic to $X$ and is endowed with this modified metric that we shall denote $g^{\delta,\varepsilon}$. We remark that $\Vol_{g^{\delta,\varepsilon}}(X)<\Vol_{g_0}(X)$ and that the identity map provide a Gromov-Hausdorff $\varepsilon'$-approximation for $\varepsilon'<\varepsilon_1(n,i_0,\kappa)$. Now observe that if we denote by $\lambda(\delta)$ the lowest eigenvalue for the Dirichlet problem on $S^n(f(\varepsilon))\setminus B(z_0,\arcsin(\delta))$ we have that:
\begin{equation}\label{SpC-mushroom1}
\lim_{\delta\f0}\lambda(\delta)=0
\end{equation}
One can be more precise: for $\delta\f0$ we have that
$$\left(\int_0^{\pi\,f(\varepsilon)}\left[\sin\left(\frac{t}{f(\varepsilon)}\right)\,f(\varepsilon)\right]^{n-1}dt\right)\lambda(\delta)\sim$$
$$\sim\left\{
\begin{array}{c}
[-\log(\arcsin(\delta))]^{-1},\;\mbox{ se }n=2;\\
(n-2)\cdot[\arcsin(\delta)]^{n-1},\;\mbox{ se } n\geq 3.
\end{array}\right.
$$
for details we refer to \cite{Cha}, Chapter II, \S5, Theorem 6. Since we can extend the corresponding eigenfunction to an eigenfunction of the closed eigenvalue problem on $(X,g^{\delta,\varepsilon})$ this shows that for $\delta\f0$ we have a family of Riemannian manifolds, satisfying the assumptions of Theorem \ref{SpC-MT2}, and such that $\frac{\lambda_1(X,g^{\delta,\varepsilon})}{\lambda_1(X,g_0)}\f0$.
\end{exmp}

\section*{Appendix. Quantitative Sobolev inequalities}

The results of this appendix are due to S. Gallot. However as S. Gallot's original results were published in a short note in the \textit{Comptes Rendus de l'Acad\'emie des Sciences} (see \cite{Ga1}), the original proofs are rather dense and we found useful to give more explanations about the method and more detailed proofs.\newline\newline
Let us consider any compact Riemannian manifold $(M,g)$ (without boundary), whose volume will be denoted by $\Vol_g(M)$ or by $V$ according to the context, and whose diameter will be denoted by $\diam(M,g)$.\newline
Let us consider the Cheeger's isoperimetric constant, $h$, and the usual isoperimetric constant, $C$, defined by
\begin{equation}\label{SpC-hC}
h=\inf_{\Omega}\frac{\Vol_g(\partial\Omega)}{\Vol_g(\Omega)}\,,\quad C=\inf_{\Omega}\frac{\Vol_g(\partial\Omega)}{\Vol_g(\Omega)^{\frac{n-1}{n}}}
\end{equation}
where $\Omega$ runs over all domains in $M$ (with piecewise regular boundary) whose volume satisfies $\Vol_g(\Omega)\leq\frac{1}{2}\Vol_g(M)$\footnote{This restriction is necessary, because otherwise the infima of $\frac{\Vol_g(\partial\Omega)}{\Vol_g(\Omega)}$ and of $\frac{\Vol_g(\partial\Omega)}{\Vol_g(\Omega)^{\frac{n-1}{n}}}$ are zero (just make the choice of $\Omega=M\setminus B(x_0,\varepsilon)$ and let $\varepsilon\f 0$).}.\newline\newline
In the euclidean space $(\R^n,\can)$, the isoperimetric constant is
$$C_*=\frac{\Vol_{\can} (\sphere^{n-1})}{(\Vol_{\can}(\Ball^n(1)))^{\frac{n-1}{n}}}$$

In the sequel we shall define the $L^p$-norms on $(M,g)$ by $$\pa u\pa_p=\left(\frac{1}{\Vol_g(M)}\,\int_M|u|^p\,dv_g\right)^{1/p}$$ and the space $H_1^p(M,g)$ ($p\geq 1$) as the completion of $C^\infty(M)$ with respect to the norm $\pa f\pa_{H_1^p}=\pa f\pa_p+\pa\nabla f\pa_p$.
\begin{lemma}\label{SpC-A} 
Let $(M,g)$ be a compact Riemannian manifold of dimension $n$.
For every domain $\Omega$ in $M$ such that $\Vol_g(\Omega)\leq\frac{1}{2}\Vol_g(M)$, and for any regular function $f\geq 0$ over $\Omega$, such that $f\mid_{\partial\Omega}=0$ we have:
\begin{itemize}
\item[(i)] $\int_\Omega|\nabla f|\,dv_g\geq C\cdot\left(\int_{\Omega}f^{\frac{n}{n-1}}\,dv_g\right)^{\frac{n-1}{n}}$;
\item[(ii)] $\left(\int_\Omega|\nabla f|^2\,dv_g\right)^{\frac{1}{2}}\geq\frac{n-2}{2(n-1)}\cdot C\cdot\left(\int_\Omega f^{\frac{2n}{n-2}}\,dv_g\right)^{\frac{n-2}{2n}}$.
\end{itemize}
\end{lemma}

\textbf{Proof.} We start proving (i). Let $\Omega_t=\{x\in\Omega|\,f(x)>t\}$ and let $\Omega_t^*$ be the open ball in $\R^n$ (centered at the origin), whose radius is determined by $\Vol_g(\Omega_t)=\Vol_{\can}(\Omega_t^*)$. We denote by $\Omega_*$ the open ball in $\R^n$ (centered at the origin) such that $\Vol_g(\Omega)=\Vol_{\can}(\Omega_*)$. We will define $A(t)=\Vol_g(\Omega_t)$ and $A^*(t)=\Vol_{\can}(\Omega_t^*)$, where $\Omega_t^*$ is the euclidean ball (centered at the origin) such that $\Vol(\Omega_t^*)=\Vol(\Omega_t)$. We construct the function $f^*:\Omega_*\f\R$ such that
$$f^*(x)=\left\{\begin{array}{c}
t\quad\quad\mbox{ when }x\in\partial\Omega_t^*;\\
\in(t-\varepsilon, t]\mbox{ when }x\in\Omega_{t-\varepsilon}^*\setminus\Omega_{t}^*;
\end{array}\right.$$
(\textit{i.e.} if $\bigcap_{\varepsilon>0}\Omega_{t-\varepsilon}^*\setminus\overline{\Omega_t^*}\neq\varnothing$, the function $f^*$ is constant and equal to $t$ on this set). It is a classical result of the symetrization method (see for instance \cite{Ber}) that
$$\int_\Omega f^{\frac{n}{n-1}}\,dv_g=\int_{\Omega_*}(f^*)^{\frac{n}{n-1}}dv_{\can}\;.$$

On the other hand, using the coarea formula (\cite{BZ}, Theorem 13.4.2) we obtain:
$$\int_\Omega|\nabla f|\,dv_g=\int_0^{\sup(f)}\Vol_g(\{f=t\})\,dt=\int_0^{\sup(f)}\Vol_g(\partial\Omega_t)\,dt\geq$$
$$\geq\int_0^{\sup(f)} C\cdot A(t)^{\frac{n-1}{n}}\,dt=\frac{C}{C_*}\int_0^{\sup(f^*)}C_*A^*(t)^{\frac{n-1}{n}}\,dt=$$$$=\frac{C}{C_*}\int_0^{\sup(f^*)}\Vol_{\can}(\partial\Omega_t^*)\,dt$$
where the last equality comes from the fact that we are in the equality-case for the isoperimetric inequality in $\R^n$ and where we intend $\int_0^{\sup(f)}$ as the integral on the set $[0,\sup(f)]\setminus\mathcal S_f$ where $\mathcal S_f$ is the set of singular values of $f$ which has measure zero by Sard's theorem. It follows that,
$$\int_\Omega|\nabla f|\,dv_g\geq\frac{C}{C_*}\int_0^{\sup(f^*)}\Vol_{\can}(\{f^*=t\})\,dt=\frac{C}{C_*}\int_{\Omega_*}|\nabla f^*|\,dv_{\can}$$
because the symmetrization method certify that $f^*$ is  Lipschitz, and thus the coarea formula (\cite{BZ}, Theorem 13.4.2) applies to $f^*$.
We get
$$\int_\Omega|\nabla f|\,dv_g\geq\frac{C}{C_*}\int_{\Omega_*} |\nabla f^*|\,dv_{\can}\geq C\cdot\left(\int_{\Omega_*}(f^*)^{\frac{n}{n-1}}\,dv_{\can}\right)^{\frac{n-1}{n}}=$$
$$= C\cdot\left(\int_{\Omega} f^{\frac{n}{n-1}}\,dv_g\right)^{\frac{n-1}{n}}$$
This ends the proof of (i).\\
Next we prove (ii). Let $h=f^{\frac{2(n-1)}{n-2}}$; since $f$ is a regular function and $f\geq 0$, and since $x\f x^{\frac{2(n-1)}{(n-2)}}$ is Lipschitz on $[0,\sup(f)]$, the function $h$ is  Lipschitz on $(M,g)$ (with bounded Lipschitz constant), so it is \textit{a.e.}-differentiable, and $h\in H_1^1(M,g)$. We have:
$$|\nabla h|=\frac{2(n-1)}{n-2}\cdot f^{\frac{n}{n-2}}\cdot|\nabla f|$$
By Lemma \ref{SpC-A} (i) we know that
$$C\cdot\left(\int_\Omega h^{\frac{n}{n-1}}\,dv_g\right)^{\frac{n-1}{n}}\leq\int_\Omega|\nabla h|\,dv_{g}$$
so that,
$$C\cdot\left(\int_\Omega f^{\frac{2n}{n-2}}\,dv_g\right)^{\frac{n-1}{n}}\leq\frac{2(n-1)}{n-2}\int_\Omega f^{\frac{n}{n-2}}|\nabla f|\,dv_g$$
which implies,
$$\frac{n-2}{2(n-1)}\cdot C\cdot\left(\int_\Omega f^{\frac{2n}{n-2}}\,dv_g\right)^{\frac{n-2}{2n}}\leq\left(\int_\Omega|\nabla f|^2\,dv_g\right)^{\frac{1}{2}}$$
and this proves Lemma \ref{SpC-A} (ii). $\Box$\\

\begin{lemma}\label{SpC-B}
Let $(M,g)$ be a compact Riemannian manifold of dimension $n$.
Let $f\in C^\infty(M)$ and let $a\in\R$ be such that $\Omega_a^+=\{f> a\}$ and $\Omega_a^-=\{f< a\}$ have volume less or equal to $\frac{\Vol_g(M)}{2}$. Then if $V=\Vol_g(M)$ we have:
\begin{itemize}
\item[(i)]$\left(\frac{1}{V}\int_M|f-a|^{\frac{2n}{n-2}}\, dv_g\right)^{\frac{n-2}{2n}}\leq\frac{2(n-1)}{(n-2)\,C\, V^{-\frac{1}{n}}}\cdot\left(\frac{1}{V}\int_M|\nabla f|^2\,dv_g\right)^{\frac{1}{2}}$
\item[(ii)] Moreover if $\int_M f\,dv_g=0$, we have, $a\leq\frac{1}{V}\int_M|f|\,dv_g$, and\newline
$\left(\frac{1}{V}\int_M|f|^{\frac{2n}{n-2}}\,dv_g\right)^{\frac{n-2}{2n}}\leq\left[\frac{2(n-1)}{(n-2)\,C\,V^{-\frac{1}{n}}}+\frac{2}{h}\right]\cdot\left(\frac{1}{V}\int_M|\nabla f|^2\,dv_g\right)^{\frac{1}{2}}$
\end{itemize}
\end{lemma}

\textbf{Proof.} First we prove (i). Applying the Lemma \ref{SpC-A} (ii) to $|f-a|$ defined on $\Omega_a^+$ (resp. $\Omega_a^-$), we obtain 
$$\left(\int_M|f-a|^{\frac{2n}{n-2}}\,dv_g\right)^{\frac{n-2}{n}}\leq$$$$\leq\left(\int_{\Omega_a^+}(f-a)^{\frac{2n}{n-2}}\,dv_g\right)^{\frac{n-2}{n}}+\left(\int_{\Omega_a^-}|f-a|^{\frac{2n}{n-2}}\,dv_g\right)^{\frac{n-2}{n}}\leq$$
$$\leq\left(\frac{2(n-1)}{(n-2) C}\right)^2\cdot\int_M|\nabla f|^2\,dv_g$$
$$\left(\int_M|f-a|^{\frac{2n}{n-2}}\,dv_g\right)^{\frac{n-2}{n}}\leq\left(\frac{2(n-1)}{(n-2)C}\right)^2\cdot\int_M|\nabla f|^2\,dv_g$$
and this proves (i). Changing eventually $f$ in $(-f)$ we can suppose that $a\geq 0$. Let us now remark that, as $\int_M f\,dv_g=0$ and $\int_{M\setminus\Omega_a^{-}} (f-a)\,dv_g\geq 0$, we have
$$a\Vol_g(M)\leq 2a\,\Vol_g(M\setminus\Omega_a^-)= 2\,\int_{M\setminus\Omega_a^-} f\,dv_g-2\,\int_{M\setminus\Omega_a^-}(f-a)\,dv_g$$
From this, from the triangle inequality and from (i) we deduce, when $\int_Mf\,dv_g=0$,
$$\left(\frac{1}{V}\int_M|f|^{\frac{2n}{n-2}}\,dv_g\right)^{\frac{n-2}{2n}}\leq\left(\frac{1}{V}\int_M|f-a|^{\frac{2n}{n-2}}\,dv_g\right)^{\frac{n-2}{2n}}+a\leq$$
$$\leq\left(\frac{2(n-1)}{(n-2)CV^{-\frac{1}{n}}}\right)\cdot\left(\frac{1}{V}\int_M|\nabla f|^2\,dv_g\right)^{\frac{1}{2}}+\left(\frac{1}{V}\int_Mf^2\,dv_g\right)^{\frac{1}{2}}\leq$$
$$\leq\left(\frac{2(n-1)}{(n-2)C V^{-\frac{1}{n}}}+\frac{2}{h}\right)\cdot\left(\frac{1}{V}\int_M|\nabla f|^2\,dv_g\right)^{\frac{1}{2}}$$
where the last inequality comes from the inequality $\frac{\int_M|\nabla f|^2\,dv_g}{\int_M f^2\,dv_g}\geq\lambda_1(M,g)\geq\frac{h^2}{4}$ (here $\lambda_1(M,g)$ stands for the first nonzero eigenvalue of the Laplace-Beltrami operator of $(M,g)$) proved by J. Cheeger in \cite{Ch1}. This ends the proof of (ii). $\Box$

\begin{prop}[Sobolev inequality, \cite{Ga1}]\label{SpC-GSobolev}
Let $(M,g)$ be a compact Riemannian manifold of dimension $n$, such that $\diam(M,g)\leq D$,  $r_{\min}\cdot D^2\geq-(n-1)\alpha^2$. For every function $f:M\f\R$ in $H_1^2(M,g)$ we have:
\begin{equation*}\mathrm{(i)}\quad\left(\frac{1}{V}\int_M|f-\overline f|^{\frac{2n}{n-2}}\;dv_g\right)^{\frac{n-2}{2n}}\leq
\end{equation*}
$$\leq\left[\frac{2(n-1)}{(n-2)\Gamma(\alpha)}+\frac{2}{H(\alpha)}\right] \cdot D\cdot\left(\frac{1}{V}\int_M|\nabla f|^2\,dv_g\right)^{\frac{1}{2}}$$
\begin{equation*}
\mathrm{(ii)}\quad\left(\frac{1}{V}\int_M f^{\frac{2n}{n-2}}\,dv_g\right)^{\frac{n-2}{2n}}\leq
\end{equation*}
$$\leq\left(\frac{2(n-1)}{(n-2)\Gamma(\alpha)}+\frac{2}{H(\alpha)}\right)\cdot D\cdot\left(\frac{1}{V}\int_M|\nabla f|^2\,dv_g\right)^{\frac{1}{2}}+\left(\frac{1}{V}\int_M f^2\,dv_g\right)^{\frac{1}{2}}$$
where we denote by $\bar f$ the mean value of $f$, \textit{i.e.} $\bar f=\frac{1}{V}\int_M f\,dv_g$, where  $H(\alpha)=\alpha\left(\int_0^{\alpha/2}(\cosh(t))^{n-1}dt\right)^{-1}$ and  where 
$$\Gamma(\alpha)=\alpha\left(\int_0^\alpha\left(\frac{\alpha}{H(\alpha)}\cosh(t)+\frac{1}{n}\sinh(t)\right)^{n-1}\,dt\right)^{-\frac{1}{n}}$$ We will use the notation $B(\alpha)$ to refer to the quantity $\left(\frac{2(n-1)}{(n-2)\Gamma(\alpha)}+\frac{2}{H(\alpha)}\right)$.
\end{prop}

\begin{rmk}
Observe that the Sobolev constant $$\Sob(g)=\inf_{f\in H_1^2(M,g),\,f\,not\,const}\frac{\pa\nabla f\pa_2}{\pa f-\bar f\pa_{\frac{2n}{n-2}}}$$
is not invariant by homotheties. It is the reason why Proposition \ref{SpC-GSobolev} bounds from below the invariant quantity $\Sob(g)\cdot\diam(M,g)$ in terms of the parameter $\alpha$.
\end{rmk}

\textbf{Proof.}  The triangle inequality
$$\pa f\pa_{\frac{2n}{n-2}}\leq\pa f-\overline f\pa_{\frac{2n}{n-2}}+|\overline f|\leq\pa f-\overline f\pa_{\frac{2n}{n-2}}+\pa f\pa_2$$
easily shows that (i) $\Rightarrow$ (ii). Hence it is sufficient to prove (i).  We admit the following isoperimetric inequalities proved by S. Gallot (\cite{Ga1},\,\cite{Ga3}), valid for every compact manifold $(M,g)$ such that $r_{\min}\cdot\diam(M,g)^2\geq-(n-1)\alpha^2$ and for any domain $\Omega$ with regular boundary and volume at most $\frac{\Vol_g(M)}{2}$:
\begin{equation}\label{SpC-IC1}
\frac{\Vol_g(\partial\Omega)}{\Vol_g(\Omega)}\geq\frac{H(\alpha)}{D},\quad
\frac{\Vol_g(\partial\Omega)}{(\Vol_g(\Omega))^{\frac{n-1}{n}}\Vol_g(M)^{\frac{1}{n}}}\geq\frac{\Gamma(\alpha)}{D}
\end{equation}
hence, passing to the infimum with respect to $\Omega$:
\begin{equation}\label{SpC-IC2}
h\geq\frac{H(\alpha)}{D},\quad C\,V^{-\frac{1}{n}}\geq\frac{\Gamma(\alpha)}{D}
\end{equation}
Now we apply Lemma \ref{SpC-B} (ii) to the function $(f-\overline f)$ and we obtain:
$$\left(\frac{1}{V}\int_M|f-\overline f|^{\frac{2n}{n-2}}\,dv_g\right)^{\frac{n-2}{2n}}\leq\left[\frac{2(n-1)}{(n-2)CV^{-\frac{1}{n}}}+\frac{2}{h}\right]\cdot\left(\frac{1}{V}\int_M|\nabla f|^2\,dv_g\right)^{\frac{1}{2}}\leq$$
$$\leq\left(\frac{2(n-1)}{(n-2)\Gamma(\alpha)}+\frac{2}{H(\alpha)}\right)\cdot D\cdot\left(\frac{1}{V}\int_M|\nabla f|^2\,dv_g\right)^{\frac{1}{2}}$$
where the second inequality comes from the inequalities (\ref{SpC-IC2}). $\Box$

\section*{Appendix A. Stability of the geometric-arithmetic inequality.}

\begin{propA1}
Let $A$ be a symmetric, non-negative matrix with real entries which satisfies the conditions:
$$\det(A)\geq(1-\sqrt\eta)^2,\quad \tr(A)\leq n(1+\eta)^{\frac{2}{n}}$$
where $0<\eta\leq\frac{1}{4}$. Then,
$$\pa A-\Id\pa^2\leq 4(n-1)^2\sqrt\eta\cdot\left(1+\frac{n+10}{n}\sqrt\eta\right)$$
\end{propA1}

Before giving the proof of the previous Proposition we state and prove the following Lemmas:

\begin{lemmaA}
Let $(x_1,..,x_n)\in\R^n$ be such that $-1<x_1\leq x_2\leq\cdots\leq x_n$ and such that $\sum_1^n x_i=0$, then 
$$\prod_{i=1}^n(1+x_i)\leq 1-\frac{n}{2(n-1)}x_1^2\leq1-\frac{\sum_1^n x_i^2}{2(n-1)^2}$$
\end{lemmaA}

\textbf{Proof of Lemma A.2.} By assumption we have $x_1\leq 0$ and 
$\sum_{i=2}^n x_i=\,=-x_1=|x_1|;$
so the geometric-arithmetic inequality gives:
$$\prod_{i=2}^n(1+x_i)\leq\left[\frac{1}{n-1}\sum_{i=2}^n(1+x_i)\right]^{n-1}=\left(1+\frac{|x_1|}{n-1}\right)^{n-1}$$
Hence we obtain the inequality:
\begin{equation}\label{SpC-appA1}
\prod_{i=1}^n(1+x_i)\leq\left(1+\frac{|x_1|}{n-1}\right)^{n-1}(1-|x_1|)
\end{equation}
Computing its derivative, it comes that the function $\frac{(1-x)^n\cdot(1-(n-1)\,x)}{\left(1-\frac{n\,(n-1)}{2}\,x^2\right)}$ is bounded above by $1$ when $x\in[0,\frac{1}{n-1})$ and
 we obtain:
$$\left(1+\frac{|x_1|}{n-1}\right)^{n-1}(1-|x_1|)\leq1-\frac{n(n-1)}{2}\left(\frac{|x_1|}{n-1}\right)^2$$
and putting this estimate in (\ref{SpC-appA1}) we end the proof of the first inequality. To prove the second inequality let $x_1=-a$; the problem is to find the maximum of $h(x)=\sum_{i=1}^n x_i^2$ over the set $$D=\left\{(x_1,..,x_n)\in\R^n|\mbox{ }x_1=-a,\,\forall i\; x_i\geq -a,\;\sum x_i=0\right\}$$
As $D$ is a $(n-2)$-dimensional simplex and $h$ is convex, it attains its maximum at one vertex of $D$, \textit{i.e.} when all the $x_i$'s are equal to $a$ except a single one. This proves the second inequality because it implies that
$$\sum x_i^2\leq\max_{y\in D} h(y)=(n-1)a^2+(n-1)^2a^2=n(n-1)x_1^2\quad\Box$$

\begin{lemmaA}[Stability of the function $A\f\frac{\det(A)}{(\frac{1}{n}\tr(A))^n}$ near its maximum]
For any real, symmetric, non-negative  $(n\times n)$-matrix $A$ we have $\forall\eta'\in(0,1)$
$$1-\eta'\leq\frac{\det(A)}{\left(\frac{1}{n}\tr(A)\right)^n}\,\Rightarrow\,\pa A-\frac{1}{n}\tr(A)\cdot\Id\pa^2\leq 2(n-1)^2\eta'\left(\frac{1}{n}\tr(A)\right)^2$$
\end{lemmaA}

\textbf{Proof of Lemma A.3.} Since $\det(A)\neq 0$ all the eigenvalues  $\lambda_i$ of $A$ are strictly positive. Let $\bar\lambda=\frac{1}{n}\tr(A)=\frac{1}{n}\sum_1^n\lambda_i$.  As, by assumption, $0<\lambda_1\leq\lambda_2\leq\cdots\leq\lambda_n$, one has
$$-1\leq\frac{\lambda_1-\bar\lambda}{\bar\lambda}\leq\frac{\lambda_2-\bar\lambda}{\bar\lambda}\leq\cdots\leq\frac{\lambda_n-\bar\lambda}{\bar\lambda}$$
we can thus apply the Lemma A.2 which gives
$$(1-\eta')\leq\prod_{i=1}^n\left(1+\frac{\lambda_i-\bar\lambda}{\bar\lambda}\right)\leq1-\frac{1}{2(n-1)^2}\sum_{i=1}^n\frac{(\lambda_i-\bar\lambda)^2}{\bar\lambda^2}$$
and thus $\pa A-\bar\lambda\,\Id\pa^2=\sum_{i=1}^n(\lambda_i-\bar\lambda)^2\leq 2\,(n-1)^2\,\eta'\,\bar\lambda^2$
which proves the Lemma. $\Box$\newline\newline

\textbf {Proof of Proposition A.1.} Let $\bar\lambda=\frac{1}{n}\tr(A)$; by the assumptions and the geometric-arithmetic inequality we have:
\begin{equation}\label{SpC-appA2}
(1-\sqrt\eta)^{\frac{2}{n}}\leq\left[\det(A)\right]^{\frac{1}{n}}\leq\bar\lambda\leq(1+\eta)^{\frac{2}{n}}
\end{equation}
so that
\begin{equation}\label{SpC-appA3}
\pa(\bar\lambda-1)\Id\pa^2=n(\bar\lambda-1)^2\leq n\left[1-(1-\sqrt\eta)^{\frac{2}{n}}\right]^2\leq\frac{4}{n}\,\eta\,(1+2\sqrt\eta)^2\leq\frac{16}{n}\eta
\end{equation}
when $\eta\leq\frac{1}{4}$. By the assumptions of the Proposition A.1 we get the following inequalities:
$$\frac{\det(A)}{\left(\frac{1}{n}\tr(A)\right)^n}\geq\frac{(1-\sqrt\eta)^2}{(1+\eta)^2}=1-\eta'\;,$$
where $\eta'=\sqrt\eta\cdot\left[\frac{2-\sqrt\eta+\eta}{(1+\eta)^2}\right]\cdot(1+\sqrt\eta)$; by Lemma A.3 we deduce that:
$$\pa A-\bar\lambda\Id\pa^2\leq 2(n-1)^2\,\eta'\,\bar\lambda^2\leq 4(n-1)^2\sqrt\eta\,(1+\sqrt\eta)\,\bar\lambda^2$$
Using equation (\ref{SpC-appA2}) we see that, when $\eta\leq\frac{1}{4}$
\begin{equation*}
\pa A-\bar\lambda\Id\pa^2\leq 4(n-1)^2\sqrt\eta\,(1+\sqrt\eta)\,(1+\eta)^{\frac{4}{n}}\leq
\end{equation*}
\begin{equation}\label{SpC-appA4}
\leq 4(n-1)^2\sqrt\eta\,\left[1+\left(\frac{n+2}{n}\right)^2\sqrt\eta\right]
\end{equation}
Since $\Id$ and $A-\bar\lambda\Id$ are orthogonal, using the estimates (\ref{SpC-appA3}), (\ref{SpC-appA4}), we obtain:
$$\pa A-\Id\pa^2=\pa A-\bar\lambda\Id\pa^2+\pa\bar\lambda\Id-\Id\pa^2\leq 4(n-1)^2\sqrt\eta\left(1+\frac{n+10}{n}\sqrt\eta\right).\mbox{ }\Box$$ 

\section*{Appendix B. Estimates  for the function $\xi$. }

\begin{lemmaB}
The infinite product $\prod_{i=0}^\infty\left(1+\frac{\beta^i}{\sqrt{2\beta^i-1}}\,x\right)^{\beta^{-i}}$ is converging for every $x\in\R^+$ and $\beta=\frac{n}{n-2}$, to a continuous function $\xi(x)$ which satisfies:
$$\xi(x)\leq e^{\frac{n}{2}\left(\frac{x}{1+x}\right)}\,(1+x)^{\frac{n}{2}},\quad\forall x\geq 0.$$
\end{lemmaB}

\textbf{Proof.} We apply the equality $(1+ax)=(1+x)\left(1+(a-1)\frac{x}{x+1}\right)$, which gives the following estimate:
\begin{equation}\label{SpC-appB1}
(1+ax)\leq(1+x)\cdot e^{(a-1)\,\frac{x}{1+x}}
\end{equation}
If we apply the estimate (\ref{SpC-appB1}) for $a=\beta^{\frac{i}{2}}$ we obtain
$$1+\frac{\beta^i}{\sqrt{2\beta^i-1}}x\leq(1+\beta^{\frac{i}{2}}x)\leq(1+x)e^{\frac{x}{1+x}\,(\beta^{\frac{i}{2}}-1)}$$
so we get for the infinite product:
$$\prod_{i=0}^\infty\left(1+\frac{\beta^i}{\sqrt{2\beta^i-1}} x\right)^{\beta^{-i}}\leq	\left(\prod_{i=0}^\infty\left[(1+x)\,e^{-\frac{x}{1+x}}\right]^{\beta^{-i}}\right)\cdot\left(\prod_{i=0}^\infty e^{\frac{x}{1+x}\beta^{-\frac{i}{2}}}\right)=$$
$$=\left[(1+x) e^{-\frac{x}{1+x}}\right]^{\sum_{i=0}^\infty\frac{1}{\beta^i}}\cdot\exp\left(\frac{x}{1+x}\cdot\sum_{i=0}^\infty\frac{1}{\beta^{\frac{i}{2}}}\right)$$
since $\sum_0^\infty\frac{1}{\beta^i}=\frac{1}{1-\frac{1}{\beta}}=\frac{n}{2}$ (since we have chosen $\beta=\frac{n}{n-2}$) and 
$$\sum_{i=0}^\infty\frac{1}{\beta^{\frac{i}{2}}}=\frac{1}{1-\frac{1}{\sqrt\beta}}=\frac{1+\frac{1}{\sqrt\beta}}{1-\frac{1}{\beta}}=\frac{n}{2}\left(1+\sqrt\frac{n-2}{n}\right)\leq n$$
we deduce that $\xi(x)\leq e^{\frac{n}{2}\,\left(\frac{x}{1+x}\right)}(1+x)^{\frac{n}{2}}$. $\Box$\newline

\begin{lemmaB}
\begin{itemize}
\item[(i)] for every $x\in\R^+$, $\xi(x)\leq(1+x)^n$;
\item[(ii)] for every $x\in[1,+\infty)$, $\xi(x)\leq(4\,e)^{\frac{n}{4}}\,x^{\frac{n}{2}}$.
\end{itemize}
\end{lemmaB}

\textbf{Proof of (i).} Bounding from above the derivative of the exponential function, one obtains, for every $t\in\R^+$, $e^t-e^0\leq t\,e^t$ and thus:
$$e^{\frac{x}{1+x}}-1\leq\frac{1}{1+x}\,e^{\frac{x}{1+x}},$$
which leads to $e^{\frac{x}{1+x}}\leq(1+x)$. From this and from the Lemma B.1, we deduce that
$$\xi(x)\leq e^{\frac{n}{2}\,\left(\frac{x}{1+x}\right)}(1+x)^{\frac{n}{2}}\leq(1+x)^{n}.\quad\Box$$
\textbf{Proof of (ii).} As $t\f\frac{e^t}{t}$ is decreasing on $[\frac{1}{2},1]$ it comes that 
$$\frac{e^t}{t}\leq\frac{e^{\frac{1}{2}}}{\frac{1}{2}}=2\sqrt e\,.$$
When $x\geq 1$, then $\frac{x}{1+x}\in[\frac{1}{2},1)$ and thus
$$e^{\frac{x}{1+x}}\,\left(\frac{1+x}{x}\right)\leq2\sqrt e.$$
A direct consequence is the estimate:
$$\xi(x)\leq e^{\frac{n}{2}\,\left(\frac{x}{1+x}\right)}(1+x)^{\frac{n}{2}}=\left[e^{\left(\frac{x}{1+x}\right)}\left(\frac{1+x}{x}\right)\right]^{\frac{n}{2}}\,x^{\frac{n}{2}}\leq(2\sqrt e)^{\frac{n}{2}}\,x^{\frac{n}{2}}.\quad\Box$$

\section*{Appendix C. A counterexample.}

\begin{lemmaC}
In general it is not possible to find a bound of the type:
$$\sup_{f\in\mathcal A(\lambda)\setminus\{0\}}\left(\frac{\pa\Delta f\pa_p}{\pa f\pa_p}\right)\leq\lambda.$$
\end{lemmaC}

\textbf{Proof.}  Let us consider, for instance, the case where $(M,g)=(\mathbb S^n,\can)$. We will denote by $\Delta_{\can}$ the corresponding Laplace-Beltrami operator. Let $\lambda=2(n+1)$, then $\mathcal A(\lambda)$ is the direct sum of the eigenspaces corresponding to of first three eigenvalues of $\Delta_{\can}$, \textit{i.e.} $\mathcal A(\lambda)=E_0\oplus E_1\oplus E_2$ where:
\begin{itemize}
\item $E_0$ is the set of the constant functions, it is the eigenspace relative to the eigenvalue $\lambda_0=0$ and $\dim(E_0)=1$;
\item $E_1$ is the space generated  by $f_1,..,f_{n+1}$, where $f_i(x)=x_i$ (here $x$ are the cartesian coordinates for $\mathbb S^n$). $E_1$ is the eigenspace corresponding to $\lambda_1=n$ and its dimension is $\dim(E_1)=n+1$;
\item $E_2$ is the space generated by the functions of the form:
$$f:\mathbb S^n\f\R,\quad f(x)=Q(x)$$
where $Q$ is a quadratic form with trace equal to zero. $E_2$ is the eigenspace corresponding to $\lambda_2=2(n+1)$ and $$\dim(E_2)=\frac{(n+2)(n+1)}{2}-1$$ (\textit{i.e.} the dimension of the $(n+1)\times(n+1)$ symmetric matrices with trace equal to zero).
\end{itemize}
We can consider the functions:
$$\varphi_0(x)=-\left(\frac{n-1}{2(n+1)}\right),\mbox{ }\varphi_0\in  E_0$$
$$\varphi_2(x)=\frac{1}{n+1}\left[nx_1^2-x_2^2-...-x_{n+1}^2\right]=\frac{1}{(n+1)}\left[(n+1)x_1^2-\sum_{i=1}^{n+1}x_i^2\right]$$
(hence, since we are restricted to $\mathbb S^n$, $\varphi_2(x)=x_1^2-\frac{1}{(n+1)}$). We remark that $\varphi_2\in E_2$ and that $\pa\varphi_2\pa_\infty=\frac{n}{n+1}$. We define $u=\varphi_0+\varphi_2$; then $u\in\mathcal A(\lambda)$ and $u(x)=x_1^2-\frac{1}{(n+1)}-\frac{n-1}{2(n+1)}=x_1^2-\frac{1}{2}$, hence $\pa u\pa_\infty=\frac{1}{2}$. On the other hand $\Delta_{\can} u=\Delta \varphi_2=2(n+1)\varphi_2$, so
$$\pa\Delta_{\can} u\pa_\infty=2(n+1)\pa\varphi_2\pa_\infty=2n$$
so we have:
$$\frac{\pa\Delta_{\can} u\pa_\infty}{\pa u\pa_\infty}=4n>\lambda=2(n+1)$$
Since the ratios $\frac{\pa\Delta_{\can} u\pa_{2k}}{\pa u\pa_{2k}}$ converge to $\frac{\pa\Delta_{\can} u\pa_\infty}{\pa u\pa_\infty}$ when $k\f\infty$, there are infinite values of $k$ for which $\frac{\pa\Delta_{\can} u\pa_{2k}}{\pa u\pa_{2k}}>\lambda$. $\Box$\newline

\end{document}